\documentclass[10pt]{amsart}
\usepackage[utf8]{inputenc}
\usepackage[a4paper,left=3cm, right= 3cm, top=3cm, bottom= 3cm]{geometry}
\usepackage{amsmath}
\usepackage{amsfonts}
\usepackage{amssymb}
\usepackage{amsthm}
\usepackage{bbold}
\usepackage{color}
\usepackage{enumerate}
\usepackage{amsmath, amssymb, amstext,amsthm}
\usepackage[toc,page]{appendix}
\usepackage{color}
\usepackage{float}

\usepackage{xcolor}
\usepackage{amsaddr}

\usepackage[active]{srcltx}
\usepackage{ifpdf}
\ifpdf 
    \usepackage[pdftex,     
            plainpages=false,   
            breaklinks=true,    
            colorlinks=true,
            pdftitle=My Document
            pdfauthor=My Good Self
           ]{hyperref} 
\else 
    \usepackage{hyperref}       
\fi

\DeclareMathOperator{\dv}{div}
\newcommand{\tor}{\mathbb{T}^2}

\newcommand{\antidvp}{\dv^{-1}_{\mathbb{T}^2}}
\newcommand{\antidvx}{\dv^{-1}_X}
\newcommand{\px}{\mathbb{P}_X}

\DeclareMathOperator{\supp}{supp}

\theoremstyle{plain}
\newtheorem{thm}{Theorem}[section]
\newtheorem{lemma}[thm]{Lemma}
\newtheorem{prop}[thm]{Proposition}

\newtheorem{cor}[thm]{Corollary}

\theoremstyle{definition}
\newtheorem{dfn}[thm]{Definition}
\newtheorem{rem}[thm]{Remark}

\newcommand{\dx}{\mathrm{d}x}

\newcommand{\dz}{\mathrm{d}z}
\newcommand{\ds}{\mathrm{d}s}

\newcommand{\di}{\mathrm{d}x_1}
\newcommand{\dii}{\mathrm{d}x_2}

\newcommand{\curl}{\operatorname{curl}}

\newcommand{\sumk}{\sum_k}
\newcommand{\tme}{\operatorname{time}}
\newcommand{\quadr}{\operatorname{quad}}
\newcommand{\Rtr}{\overset{\circ}{R_0}}
\newcommand{\bl}{\mathcal{B}}
\newcommand{\vl}{\mathcal{V}}
\newcommand{\Rc}{\mathcal{R}}

\def\R{\mathbb{R}}

\def\N{\mathbb{N}}
\def\Z{\mathbb{Z}}
\def\P{\mathbb{P}}

\def\div{\operatorname{div}}
\def\curl{\operatorname{curl}}

\def\supp{\mathrm{supp\, }}
\def\e{\varepsilon}

\title[Compactly supported anomalous solutions for 2D Euler in Hardy spaces]{Compactly supported anomalous weak solutions for 2D Euler equations with vorticity in Hardy spaces} 

\author{Miriam Buck}
\address{Technische Universität Darmstadt, Department of Mathematics, 64285 Darmstadt, Germany\\ mbuck@mathematik.tu-darmstadt.de}
\author{Stefano Modena}
\address{Gran Sasso Science Institute, 67100 L'Aquila, Italy\\ stefano.modena@gssi.it\\ ORCID: 0000-0002-2890-4285 }

\begin{document}

\begin{abstract}
In \cite{buck2023non}, we constructed by convex integration examples of energy dissipating solutions to the 2D Euler equations on $\R^2$ with vorticity in the real Hardy space $H^p(\R^2)$. In the present paper, we develop tools that significantly improve the result in \cite{buck2023non} in two ways: Firstly, we achieve vorticities in $H^p$ in the optimal range $p\in (0,1)$ compared to $(2/3,1)$ in our previous work. Secondly, the solutions constructed here possess compact support and in particular preserve linear and angular momenta.
\end{abstract}
\maketitle

\section*{Acknowledgements}
The authors would like to thank Phil Isett for a stimulating email exchange about some aspects of the present paper.

\section{Introduction}
\subsection{Main result and literature overview}

 In this paper we consider the 2-dimensional incompressible Euler equations on the full space $\R^2$
\begin{align}\label{2D Euler}
\begin{cases}
\partial_t u + \dv(u\otimes u) + \nabla p = 0,\\
\dv u = 0,\\
u(\cdot, 0) = u_0,
\end{cases}
\end{align}
where $u: [0,1]\times\R^2\rightarrow\R^2$ is the velocity field of some fluid and $p:[0,1]\times\R^2\rightarrow\R$ is the corresponding (scalar) pressure. We prove the following theorem:
\begin{thm}[Main Theorem]
\label{thm:main}
Let $0<p<1$. For any energy profile $e\in C^\infty\left([0,1];\left[\frac{1}{2},1\right]\right)$ there exists a solution $u\in C([0,1],L^2(\R^2))$ to \eqref{2D Euler} with
\begin{enumerate}[(i)]
\item $u(t,\cdot)=0$  on $\R^2\setminus B_1(0)$ a.e. for every $t\in [0,1]$,
\item $u$ preserves linear and angular momentum, i.e.
\begin{align*}
\int_{\R^2} u(t,x) \,\dx = \int_{\R^2} u(0,x)\,\dx, \qquad \int_{\R^2}  x \times u(t,x) \,\dx = \int_{\R^2} x \times u(0,x)\,\dx \text{ for every } t\in [0,1],
\end{align*} 
\item $\int_{\R^2} |u|^2(t,x) \,\dx = e(t)$,
\item $\curl u \in C([0,1],H^p(\R^2))$.
\end{enumerate}
In particular, there exist energy dissipating solutions $u \in C_t L^2_x$ to \eqref{2D Euler} with $\curl u \in C_t H^p_x$.
Furthermore, for energy profiles $e_1,e_2$ such that $e_1=e_2$ on $[0,t_0]$ for some $t_0\in [0,1]$, there exist two solutions $u_1,u_2$ satisfying $(i)-(iv)$ with $u_1(t)=u_2(t)$ for $t\in [0,t_0]$.
\end{thm} 

From this theorem one immediately deduces the following nonuniqueness result:

\begin{cor}
Let $0<p<1$. There are two admissible (in the sense that the total kinetic energy is non-increasing in time) solutions $u_1, u_2 \in C([0,1]; L^2(\R^2))$ with $\curl u_1, \curl u_2 \in C([0,1]; H^p(\R^2))$ with the same initial datum $u_1|_{t=0} = u_2|_{t=0}$.
\end{cor}
\begin{proof}
The proof follows immediately from Theorem \ref{thm:main}, picking two non-increasing energy profiles $e_1, e_2$ which coincide on $[0,1/2]$ and are different from each other on $(1/2,1]$. 
\end{proof}
Theorem \ref{thm:main} is motivated by the following question:\\

\textbf{Q}: \noindent \emph{For $u_0\in L^2(\R^2)$ with $\curl u_0\in X$ for some function space $X$, does there exist a unique solution $u\in C([0,1],L^2(\R^2))$ to \eqref{2D Euler} with $\curl u \in C([0,1],X)$ and initial datum $u_0$?}\\

Such question originates from the observation that formally  \eqref{2D Euler} can be stated as a transport equation for the vorticity $\omega = \curl u$, i.e.
\begin{align}\label{voriticty eq}
\begin{cases}
\partial_t\omega + u\cdot \nabla\omega = 0,\\
u= \nabla^\perp\Delta^{-1}\omega.
\end{cases}
\end{align}
From \eqref{voriticty eq} one sees that for any $p\in [1,\infty]$ the $L^p$ norm of the vorticity of any smooth solution to \eqref{2D Euler} is conserved in time. Since very little is known for the case $X=L^p$, it is interesting to ask \textbf{Q} for more general function spaces $X$. \\

There have been several important results providing partial answers to the question above.
We mention a few of them; for a more detailed overview we refer to the introduction in \cite{buck2023non} and to the references therein.
The first known result addresses the case $X=L^\infty$ and dates back to Yudovich \cite{yudovich1962some,yudovich1963non}. It states that for any initial datum $u_0\in L^2$ with $\omega_0\in L^1\cap L^\infty$, there exists a unique global solution $u\in C_tL^2_x$ with $\omega\in L^\infty_t(L^1_x\cap L^\infty_x)$ to \eqref{voriticty eq}. Although, in \cite{vishik2018instability1,vishik2018instability2} Vishik gave a negative answer to {\bf Q} for $X=L^p$, $p<\infty$, proving nonuniqueness in the class of solutions having vorticity $\omega \in L^\infty_t(L^p_x)$, however not for the Euler system \eqref{2D Euler} (or \eqref{voriticty eq}), but for the Euler system \eqref{2D Euler} with a $L^1_t (L^1_x \cap L^p_x)$ external force  (thus allowing for an additional ``degree of freedom'').\\

In \cite{brue2021nonuniqueness}, Bruè and Colombo address \textbf{Q} for the case that $X$ is the Lorentz space $X=L^{1,\infty}$.  They show the existence of a sequence $\{u_n\}_n$ of approximate solutions to \eqref{2D Euler} converging strongly in $L^2$ to an anomalous weak solution $u$ to \eqref{2D Euler} (therefore providing an example of nonuniqueness) and whose corresponding vorticities $\{\curl u_n\}_n$ build a Cauchy sequence in $L^{1,\infty}$ which thus has a limit $\omega$ in $L^{1, \infty}$. However, since $L^{1,\infty}$ is not a space of distributions (i.e. $L^{1,\infty}\not\hookrightarrow\mathcal{D}'$),  it is not clear whether and in what sense the distributional vorticity of the solution $u$ (or, in other words, the distributional limit of $\curl u_n$) coincides with the $L^{1,\infty}$ limit $\omega$. In general, there is no connection between distributional limit and limit in $L^{1,\infty}$ (see, for instance, the example in Section 1.4 in \cite{buck2023non}.)\\

The proof of the result in \cite{brue2021nonuniqueness} relies on an \emph{intermittent convex integration scheme}. As a general principle in intermittent convex integration, one expects such schemes to produce  weak solutions to the $d$-dimensional Euler equations having vorticity in $L^p$ if

\begin{equation}
\label{eq_lp_vorticity}
p < \frac{2d}{d+2}.
\end{equation}
The reasoning behind this principle is the fact that for such $p$ the Sobolev embedding into $L^2$ fails, i.e. $W^{1,p}\not\hookrightarrow L^2$;  see the introduction in \cite{buck2023non} for a more detailed explanation.
In particular, in dimension $d=2$, it is not possible with the current techniques to construct solutions $u$ with $\curl u \in L^p$, not even for $p=1$. The result by Bruè and Colombo motivated us in \cite{buck2023non}  to see if convex integration methods can be used to show non-uniqueness of weak solutions to \eqref{2D Euler} with vorticity in a function space $X$ that is ``weaker'' than $L^1$ in terms of integrability, but at the same time it does embed into $\mathcal{D}'$ in order to avoid the topological issues that we mentioned for $L^{1, \infty}$.  One natural choice for such space $X$ are the real Hardy spaces $H^p$ for $p<1$, which actually do embed into $\mathcal{D}'$. In \cite{buck2023non}, we proved an analogue of Theorem \ref{thm:main}, albeit only for $p\in (2/3,1)$ and solutions $u$ having merely $L^2$ decay; see Subsection \ref{subsec: technical part} for an elaboration on the improvements in the present paper.  We collect all necessary information about real Hardy spaces in Subsection \ref{subsec: hardy spaces}. For a more detailed summary of why real Hardy spaces are the right substitute for $L^p$ spaces in that context, see the introduction in \cite{buck2023non}.

\subsection{Novelties and technical challenges}\label{subsec: technical part}

We already addressed the above question \textbf{Q} with $X=H^p(\R^2)$ in \cite{buck2023non}. There  we constructed energy dissipating solutions $u\in C_tL_x^2$ to the 2D Euler equations with vorticity $\curl u \in H^p(\R^2)$ for $2/3<p<1$. In the present paper, we resolve questions that remained left open in \cite{buck2023non}, leading to a significant improvement of the result in \cite{buck2023non}. We tackle the following issues:

\begin{enumerate}
    \item Concerning the decay at infinity, the solutions in  \cite{buck2023non} are merely $C_tL_x^2$-functions. With our method at that time, achieving solutions in $C_t(L_x^1\cap L_x^2)$ was out of reach, let alone solutions with well-defined linear and angular momentum.  We believe that with some effort, it is possible to show that solutions constructed in  \cite{buck2023non} are in fact in $C_tL_x^s$ for $1<s\leq 2$, excluding $s=1$, although it is not proven there. In the present paper, however, we construct \textit{solutions having compact support} in space, hence $u\in C_tL_x^s$ for $1\leq s \leq 2$ is automatic and $u$ has well-defined linear and angular momentum. Observe that for any weak solution $u$ of \eqref{2D Euler} linear and angular momentum are conserved in time if they are well-defined,  see for example \cite{isett2016fullspace}, Proposition 3.1.
    \item In \cite{buck2023non}, we did only consider the case $2/3<p<1$. Restricting to that range of $p$ is a great simplification: indeed, for a compactly supported function $\varphi$, $\curl\varphi\in H^p(\R^2)$ always holds for $p$ in that range. On the other hand, for $p\leq 2/3,$ the requirement $\curl\varphi\in H^p(\R^2)$ is far from trivial, even for $\varphi\in C^\infty_c$. In the present construction, we reach the \textit{optimal range $0<p<1$}. (Note that $p\geq1$ is not compatible with intermittent convex integration as pointed out in \eqref{eq_lp_vorticity}.) 
\end{enumerate}

Therefore, the present construction yields solutions with the best possible decay (i.e. compactly supported functions) and the greatest possible flexibility for $p$ (the full range $0<p<1$) in the class of $L^2$-solutions to \eqref{2D Euler} with vorticity in some space $H^p(\R^2)$. 

\section{Comments on the proof of Theorem \ref{thm:main}}
In this section, we will further comment on  the above mentioned issues and how we overcome them. First, let us spend some words about the convex integration scheme used here, thereby already introducing some notation needed for the subsequent subsections. As in \cite{buck2023non}, we use an intermittent convex integration technique in the spirit of De Lellis and Székelyhidi works on the 3D Euler equations in the framework of Onsager's Theorem (see \cite{de2009euler,de2014dissipative,de2013dissipative,isett2018proof,buckmaster2018}) and Buckmaster and Vicol works on the Navier-Stokes equations (see \cite{buckmaster2019nonuniqueness}). The outline in all of these schemes is an iterative construction where, starting from an initial approximate solution, one adds fast oscillating perturbations with a higher frequency $\lambda_n \to \infty$ with respect to the typical frequencies $\lambda_{n-1}$ in the previous approximation. In case of the Euler equation, given an approximate solution $(u_{n-1},p_{n-1},R_{n-1})$ with error term on the right hand side
\begin{align}\label{eq: error previous}
\begin{cases}
\partial_t u_{n-1} + \dv(u_{n-1}\otimes u_{n-1}) + \nabla p_{n-1}= - \dv R_{n-1},\\
\dv u_{n-1} = 0,
\end{cases}
\end{align}
one makes the Ansatz
\begin{align*}
u_{n}(t,x)&= u_{n-1}(t,x) + w_n(t,x) 
\end{align*}
with
\begin{align}\label{intro perturbations}
w_n(t,x) &= u_n^p(t,x) + u_n^t(t,x) + \text{corrector terms},\nonumber\\
u_n^p(t,x) &= a_{n-1}(t,x) W_{\mu_n}(t,\lambda_n x) \quad\text{"principal part"},\nonumber\\
u_n^t(t,x) &=\P\left(a_{n-1}^2(t,x) Y_{\mu_n}(t,\lambda_n x)\right) \quad \text{"time corrector"},
\end{align}
where $\P$ is the Leray projector, $\xi$ is a given direction and $\left\lbrace W_\mu\right\rbrace_{\mu>0}, \left\lbrace Y_\mu\right\rbrace_{\mu>0}$ are families of "building blocks" with scaling
\begin{align}\label{eq: scaling building blocks}
\|\nabla^kW_\mu\|_{L^s} \approx \mu^{\frac{d}{2}-\frac{d}{s}+k}, \quad \|\nabla^k Y_\mu\|_{L^s} \approx \mu^{d-\frac{d}{s}+k}
\end{align}
and the property
\begin{align}\label{eq: cancellation divergence}
\begin{cases}
\partial_t Y_{\mu}(t,\lambda_n x) + \div\left(W_{\mu}(t,\lambda_n x)\otimes W_{\mu_n}(t,\lambda_n x)\right) = 0,\\
\div W_\mu = 0.
\end{cases}
\end{align}
Here, $\lambda_n, \mu_n$ are the oscillation and concentration parameter, respectively. Furthermore, $a_{n-1}$ is a  slowly varying coefficient with $ a_{n-1} \approx |R_{n-1}|^{1/2}$. The role of the principal part $u_n^p$ is the cancellation of the error term $R_{n-1}$:  the interaction of $u_n^p$ (having frequencies $\lambda_n$) with itself from the nonlinearity of the equation produces a term having frequencies $\approx \lambda_{n-1}$ and it allows therefore for the cancellation of the previous error. On the other hand, the purpose of the time corrector $u_n^t$ is to cancel the error $\dv(W_{\mu_n}(t,\lambda_n\cdot)\otimes W_{\mu_n}(t,\lambda_n\cdot))$ which is of order $\lambda_n$ in $L^1$ and that cancellation is crucial for the iteration step; see also \cite{modena2020convex} that follows the same strategy in the context of the transport equation. 
The principal part $u_n^p$ is not divergence free. The role of the further corrector terms is to compensate for the divergence of $u_n^p$ and make $w_n$ divergence-free. One can do this by using that the fast oscillating part  $W_{\mu_n}(t,\lambda_n x)$ is divergence-free and has vanishing mean value and therefore can be written as a potential. For the time corrector this is not possible: the building block $Y_{\mu_n}$ that is built to cancel the quadratic divergence term is quadratic itself and therefore has non-vanishing mean value. Therefore one includes the Leray projector $\P$ in the definition of $u_n^t$.

\subsection{Obtaining compact support}
The reason why in \cite{buck2023non} we were unable to achieve solutions with compact support (or at least solutions with $L^1$ decay) was grounded in the fact that we were using the Leray projector on two occasions. One is the time corrector $u_n^t$ introduced above, the second one is another corrector that we will futher comment on below. Since $\P$ is a nonlocal operator, we could not hope for those parts of the perturbation to have compact support despite the fact that the arguments inside $\P$ are compactly supported. At the same time, there are no $L^1$-bounds for $\P$ available.

\subsubsection{The time corrector}\label{subsec: time corrector}
Recalling Definition \eqref{intro perturbations}, the quadratic building block $Y_{\mu_n}$ for the time corrector $u_n^t$ is not divergence-free nor can it be made divergence-free by adding a small corrector as for $u_n^p$. Therefore, in \cite{buck2023non} we include the Leray projector $\P$ into the definition of $u_n^t$. In the present paper we choose a different route and take advantage of the structure $a_{n-1}^2 Y_{\mu_n}(t,\lambda_n x)$, which is a product of a compactly supported function and a periodic fast oscillating function. Noting that $$\left|\int a_{n-1}^2 Y_{\mu_n}(t,\lambda_n x) \,\dx\right| \approx \frac{1}{\omega_n}\ll 1,$$ where $\omega_n$ is the phase speed, see Section \ref{sec: perturbations}, we consider $$a_{n-1}^2 \left(Y_{\mu_n}(t,\lambda_n x) - \frac{1}{\omega_n}\xi\right),$$ where $\xi\in\R^2$ is a given direction, and the second (periodic) factor has vanishing mean value. Hence there exists an inverse Laplacian on $\tor$ for the second factor and we can make the following Ansatz, where we replace $u_n^t$ by a perturbation $w_n^t$ that "contains" $u_n^t$ as we shall see below, namely
\begin{align*}
w_n = u_n^p + w_n^t
\end{align*}
with
\begin{align}\label{eq time corrector 2}
w^t_n(t,x) = \nabla^\perp\curl \left(a_{n-1}^2\Delta^{-1} \left(Y_{\mu_n}(t,\lambda_n x)-\frac{1}{\omega_n}\xi\right)\right),
\end{align}
or, more generally
\begin{align}\label{eq time corrector 3}
w^t_n(t,x) = \nabla^\perp\curl\Delta^N\curl \left(a_{n-1}^2\Delta^{-(N+1)} \left(Y_{\mu_n}(t,\lambda_n x)-\frac{1}{\omega_n}\xi\right)\right).
\end{align}
The advantage of \eqref{eq time corrector 2} or \eqref{eq time corrector 3} is that any orthogonal gradient is automatically divergence free. The importance of using a higher order Laplacian as in \eqref{eq time corrector 3} lies in the fact that we need to control higher order momenta of the perturbations to achieve $\curl w_n^t\in H^p$ and this will become apparent in Section  \ref{sec: issues hardy spaces}. Let us for the moment assume that $w_n^t$ is defined by \eqref{eq time corrector 2}, then, noting that $\nabla^\perp\curl = \Delta - \nabla\div$, this yields
\begin{align}\label{eq time corrector 4}
w^t_n(t,x) &= a_{n-1}^2Y_{\mu_n}(t,\lambda_n x) - \frac{a_{n-1}^2}{\omega_n}\xi -\nabla\div \left(a_{n-1}^2 \left(Y_{\mu_n}(t,\lambda_n x)-\frac{1}{\omega_n}\xi\right)\right)\\
&\hspace{0,3cm} + \underbrace{\Delta a_{n-1}^2 \Delta^{-1} \left(Y_{\mu_n}(t,\lambda_n x)-\frac{1}{\omega_n}\xi\right) + \nabla a_{n-1}^2\cdot
\nabla\Delta^{-1} \left(Y_{\mu_n}(t,\lambda_n x)-\frac{1}{\omega_n}\xi\right)}_{\text{ lower order terms of order }\leq\frac{1}{\lambda_n}}.\nonumber
\end{align}
The first term is exactly the time corrector $u_n^t$ (without the Leray projector) that is needed to cancel the quadratic divergence term as in \eqref{eq: cancellation divergence}, the second one is the part we added to make the fast oscillating part mean value free. The third summand in \eqref{eq time corrector 4} is a gradient which can be taken care of by the pressure term in \eqref{2D Euler}. The lower order terms all gain at least a factor $\frac{1}{\lambda_n}$. However, this is not sufficient to estimate their time derivatives in $L^1$, which is necessary for the control of the new error term. But, noting that the lower order terms are again products of a compactly supported function and a fast oscillating function, one can expand $w_n^t$ by some correctors (which arise from applying once more $\Delta^{-1}$ to the oscillating parts of the surplus expressions) that will cancel them: the refined Ansatz for $w_n^t$ is then
\begin{align*}
w^t_n(t,x) &= \nabla^\perp\curl \left(a_{n-1}^2\Delta^{-1} \left(Y_{\mu_n}(t,\lambda_n x)-\frac{1}{\omega_n}\xi\right) \right)\\
&\hspace{0,3cm}+ \nabla^\perp\curl \left( \Delta a_{n-1}^2\Delta^{-2}\left(Y_{\mu_n}(t,\lambda_n x)-\frac{1}{\omega_n}\xi\right) + \nabla a_{n-1}^2\cdot\Delta^{-1}\nabla\Delta^{-1} \left(Y_{\mu_n}(t,\lambda_n x)-\frac{1}{\omega_n}\xi\right)\right)\\
&= u_n^t - \frac{a_{n-1}^2}{\omega_n}\xi+ \nabla\div(\dots) + (\text{ lower order terms of order } \leq\frac{1}{\lambda_n^2}),
\end{align*}
i.e. this gives the desired perturbation $u_n^t(t,x)= a_{n-1}^2(t,x)Y_{\mu_n}(t,\lambda_n x)$ up to a gradient and lower order terms of order less than $\frac{1}{\lambda_n^2}$ which are again products of a compactly supported function and a fast oscillating function. Repeating such an integration-by-parts argument $M$ times defines a bilinear operator
\begin{align*}
\bl^1_M:C^\infty_c(\R^2)\times C_0^\infty(\tor)\rightarrow C^\infty_c(\R^2)
\end{align*} 
with the property that
\begin{align*}
\Delta \bl_M^1\left(a_{n-1}^2, Y_{\mu_n}(t,\lambda_n\cdot)-\frac{1}{\omega_n}\xi\right) = a_{n-1}^2 \left(Y_{\mu_n}(t,\lambda_n\cdot)-\frac{1}{\omega_n}\xi\right) +  (\text{ lower order terms of order } \leq\frac{1}{\lambda_n^M}),
\end{align*}
such that we can choose $w_n^t$ as
\begin{align}\label{eq: w time written out}
w_n^t &= \nabla^\perp\curl\left(\bl^1_M\left(a_{n-1}^2, Y_\mu(t,\lambda_n \cdot)-\frac{1}{\omega_n}\xi\right)\right)\nonumber\\
&= u_n^t- \frac{a_{n-1}^2}{\omega_n} \xi+ \nabla\div(\dots) + (\text{ lower order terms of order } \leq\frac{1}{\lambda_n^M}),
\end{align}
which is eventually enough to control all relevant error terms associated to the the time corrector. A similar argument can be carried out for the Ansatz \eqref{eq time corrector 3} with general $N\in\N$, defining a bilinear operator $\bl^{N+1}_M$. We make this precise in \mbox{Proposition \ref{prop: improved anti laplacian}.} Our final choice of $w_n^t$ will be
\begin{align}\label{eq: ut with corrector}
w_n^t &= \nabla^\perp\curl\Delta^N\left(\bl^{N+1}_M\left(a_{n-1}^2, Y_\mu(t,\lambda_n\cdot)-\frac{1}{\omega_n}\xi\right)\right)\nonumber\\
&= u_n^t + \nabla\div(\dots) + (\text{ lower order terms of order } \leq\frac{1}{\lambda_n^M})
\end{align}
with some large numbers $N, M\in\N$. The role of $N$ is to control an arbitrarily high order of the Laplacian $\Delta^N$, while the role of $M$ is to gain control over the error terms from the fast oscillation, i.e. terms of order $\lambda^{-M}$.

\subsubsection{The antidivergence corrector}
The second occasion where we use the Leray projector in \cite{buck2023non} has to do with the fact that we are working on the full space $\R^2$ rather than on a periodic domain. It is known since the paper \cite{isett2016fullspace} by Isett and Oh that one key obstruction that one faces when  carrying out convex integration techniques on the  full space is the absence of a bounded right inverse $\dv^{-1}:L^1(\R^2;\R^2)\rightarrow L^1(\R^2;\operatorname{Sym}_{2\times 2})$ for the divergence, where $\operatorname{Sym}_{2\times 2}$ is the space of symmetric matrices. Such an operator is needed to define the right hand side of \eqref{eq: error previous}. Indeed, assuming such an operator exists, then $u\in L^1(\R^2;\R^2)$, $R=\dv^{-1}u$ with $\|R\|_{L^1}\lesssim \|u\|_{L^1}$ yields
\begin{align*}
\int u \,\dx = \int \dv R\,\dx = 0, \quad \int x\times u\,\dx = \int x\times \dv R =0,
\end{align*}
i.e. $u$ needs to have vanishing linear and angular momentum.\\
 In \cite{buck2023non}, we overcome this issue by considering the underdetermined system
\begin{align*}
\partial_t u_{n-1} + \dv(u_{n-1}\otimes u_{n-1}) + \nabla p_{n-1}= - r_{n-1} -  \dv R_{n-1}
\end{align*}
instead of \eqref{eq: error previous}, which gives a higher degree of freedom. Here $r_{n-1}$ is an arbitrarily small vector field with compact support that is constructed in each iteration step alongside $R_{n-1}$; see the introduction in \cite{buck2023non} for a more detailed presentation. In order to cancel this additional error term, in \cite{buck2023non} we include in our definition of $w_n$ the corrector
\begin{align*}
v_n(t,x) = \int_0^t \P r_{n-1}(s,x)\,\ds
\end{align*}
such that $\partial_t v$ cancels $r_{n-1}$ up to a gradient. Again, the Leray projector is needed since $r_{n-1}$ is, in general, not divergence free, and therefore there is no hope to maintain the compact support or to estimate $v_n$ in $L^1$.\\
Another way to resolve this issue is presented in \cite{isett2016fullspace}. There as well the authors have to deal with the absence of a bounded right inverse $\div^{-1} : L^1 \to L^1$ with values in the symmetric matrices. They  construct an antidivergence operator which is defined only on the subset of $L^1(\R^3; \R^3)$ consisting of elements which are orthogonal to translation and rotational vector fields. In the present paper, we exploit Isett and Oh's strategy and make sure our perturbations satisfy the orthogonality condition.
Yet another possibility is the approach taken in \cite{enciso2024extension}. There, the authors construct another antidivergence operator defined  on the space of vector fields orthogonal to translations and rotations based on the potential theoretic antidivergence and an analysis of some second order symmetric differential  operator already used in \cite{de2009euler} that is designed to correct the support of the potential theoretic solution.

\subsection{Issues connected to the real Hardy spaces}\label{sec: issues hardy spaces}
In order to ensure that the solution $u$ to \eqref{2D Euler} that we produce satisfies $\curl u\in H^p(\R^2)$, we need to make sure that the vorticities of the approximate solutions $\curl u_n$ and therefore also the vorticities of the perturbations $\curl w_n$ belong to $H^p(\R^2)$, and we need an estimate showing that $(\curl u_n)_n$ is a Cauchy sequence in $H^p(\R^2)$.

\subsubsection{Accomplishing $\curl u\in H^p$}
We recall a basic property of Hardy spaces: for a function $\varphi\in C^\infty_c(\R^d)$ it holds that $\varphi\in H^p(\R^d)$ if and only if
\begin{align}\label{eq: control momenta}
\int x^\alpha \varphi\,\dx = 0
\end{align}
for every $|\alpha|\leq d(1/p -1)$. Moreover, if $\varphi\in H^p(\R^d)$ and $\supp\varphi\subset B$ for some ball $B$, it holds
\begin{align}\label{est: atoms intro}
\|\varphi\|_{H^p}\lesssim |B|^\frac{1}{p}\|\varphi\|_{L^\infty},
\end{align}
see also Section \ref{subsec: hardy spaces} for the definition of real Hardy spaces and the above estimate. For $d=2$ and $p\in (2/3,1)$, one thus has to control only the zeroth order momentum. This condition is satisfied trivially if $\varphi=\curl\psi$ for some $\psi\in C^\infty_c(\R^d)$ or even $\varphi = \curl\P\psi$.\\
 For general $p$ in the full range $(0,1)$, we notice that the vorticity of the perturbation $w_n^t$ (the part of the perturbation "containing" the time corrector $u_n^t$) defined by \eqref{eq: ut with corrector} satisfies  \eqref{eq: control momenta}  if $2N+3\geq 2(1/p -1) $. As for the time corrector above, we will replace the principal part $u_n^p$ by a perturbation $w_n^p$ carrying the principal corrector $u_n^p$ in a similar fashion: In Section \ref{sec: building blocks} we explicitly define a symmetric matrix $A=A(N)$ such that
\begin{align}\label{eq: up}
w_n^p :&= \nabla^\perp\curl\Delta^N(\div(a_{n-1}A))\\
&= u_n^p+  \text{ lower order corrector terms},\nonumber
\end{align}
and the corrector terms can be controlled in the relevant norms. Setting 
\begin{align}\label{eq: introduction of wn}
w_n = w_n^p + w_n^t
\end{align} yields $\curl w_n\in H^p(\R^2)$ if $N$ is large enough. Notice also that \eqref{eq: ut with corrector} and \eqref{eq: up} imply that $w_n^p$ and $w_n^t$ belong to the space of functions which are orthogonal to translation and rotational vector fields.

\subsubsection{Estimates in Hardy spaces via intermittency}
In order to control the quantity $\|\curl w_n^p\|_{H^p}$, we use the mechanism of \emph{concentration} or \emph{intermittency} and the notion of a Hardy space atom.
The building block matrix $A$ introduced in \eqref{eq: up} is defined via concentrated functions,
\begin{align*}
A := A_{\mu_n} :=\varphi_{\mu_n}(x) (\xi\otimes\xi^\perp + \xi^\perp\otimes\xi)
\end{align*}
where $\varphi\in C^\infty_c(\R^2)$, $\varphi_{\mu_n}$ is the periodization of the concentrated function $\mu_n^{1-(2N+3)}\varphi(\mu_n x)$ and $\xi\in\R^2$ is some given direction. In particular, $A$ and therefore also $w_n^p$ are supported on disjoint small balls of radius $\frac{1}{\mu_n}.$ The scaling is such that we keep \eqref{eq: scaling building blocks}, i.e.
$\|W_{\mu_n}\|_{L^2} \approx \|D^{2N+3}A\|_{L^2} \,\raisebox{-0.9ex}{\~{}}\, 1$.  

 Hardy space atoms are typical functions $f$ in Hardy space that have support in a ball $B$ and satisfy the cancellation property $\int_B x^\alpha f \,\dx = 0$ for every $|\alpha|\leq 2(1/p-1)$ and an $L^\infty$ estimate, see Definition \ref{dfn: hardy atoms}. Indeed, thanks to the compact support of $a_{n-1}$ and the intermittency, one can view the perturbation $w_n^p = \nabla^\perp\curl\Delta^N(\div(a_{n-1}A))$ as a linear combination of atoms, where each of them is supported on a very small ball of radius $\frac{1}{\mu_n}$, i.e.
\begin{align*}
w_n^p &= \sum \theta_j,\\
\theta_j &=  \mathbb{1}_{B_{\frac{1}{\mu_n}}(x_j)} w_n^p
\end{align*}
for some finite set of points $x_1, \dots, x_n$. The curl of each such  $\theta_j$ is a derivative of order $2N+4\geq 2(1/p-1)$ and therefore satisfies the cancellation property $\int_{B_{\frac{1}{\mu_n}}(x_j)} x^\alpha \curl \theta_j \,\dx = 0$, \mbox{$|\alpha|\leq 2(1/p-1)$.} As a result, $\curl w_n^p$ is a linear combination of atoms and thus $\curl w_n^p\in H^p$. One can use a standard estimate for atoms (see \eqref{est: atoms intro} and also Lemma \ref{lemma: hardy atoms}) on each $\curl\theta_j$, balancing $\|\curl\theta_j\|_{L^\infty}$ (estimated by \eqref{eq: scaling building blocks}) and the size of its support, which is of order $\frac{1}{\mu_n^2}$.  

\subsubsection{Estimates in Hardy spaces via a generalized triangle inequality}
For $w_n^t$, the perturbation carrying the time corrector, an estimate in $H^p(\R^2)$ is more delicate. In contrast to $w_n^p$, we cannot split $w_n^t$ into a linear combination of atoms with small supports. Although the building blocks $Y$, $W$ and the matrix $A$ share the same support (see Lemma \ref{lemma: supports} for a precise description), this support is not seen by $w_n^t$ since we are considering $Y_{\lambda_n}-\frac{1}{\omega_n}\xi$ instead of just $Y_{\lambda_n}$ and because of the nonlocal operator $\Delta^{-1}$ used in the construction of $\bl^{N+1}_M$, see Section \ref{subsec: time corrector}. One might be tempted to apply an estimate for atoms \eqref{est: atoms intro} to every summand in the second line of \eqref{eq: ut with corrector}. However, while $\curl w_n^t$ is an $H^p$-function due to the choice of $N$, that might not be true for its individual parts, and therefore, in general
\begin{align*}
\|\curl w_n^t\|_{H^p}^p \nleqslant \|\curl u_n^t\|_{H^p}^p +\|\curl\left(\frac{a_{n-1}^2}{\omega_n}\xi\right)\|_{H^p} +  \|\curl\text{(lower order terms)}\|_{H^p}^p.
\end{align*}
Observe that $\curl\nabla \equiv 0$ and the third summand in \eqref{eq: w time written out} vanishes when taking the curl.
In order to resolve that issue, we prove in Proposition \ref{Hardy for sum} an estimate for finite families of functions $(f_j)_j$ with respective supports in balls of radii $(\varepsilon_j)_j$ where $\sum_j f_j\in H^p$ without further assumptions on $(f_j)_j$, which is a proper generalisation of the estimate for atoms \eqref{est: atoms intro}. Namely,
\begin{align*}
\|\sum_{j} f_j\|_{H^p(\R^d)}^p&\leq C \sum_{j}\left(\varepsilon_j^d \|f_j\|^p_{L^\infty(\R^d)} + \max_{|\alpha|\leq d(1/p-1)}\left|\int x^\alpha f_j(x) \,\dx\right|^p\right).
\end{align*}
Note that if $f_j\in H^p(\R^d)$, the last summand vanishes and we are left with the estimate for scaled atoms as in \eqref{est: atoms intro}.

\subsection{Notation}\label{subsection: notation}

We fix some notation we are going to use in the paper.
\begin{enumerate}[(a)]
\item We denote by $e_1, e_2$ the standard basis vectors of $\R^2$.
\item For any vector $\xi=(\xi_1,\xi_2)\in\R^2$, we will denote by $\xi^\perp$ the orthogonal vector $\xi^\perp = (\xi_2, - \xi_1)$.
\item We denote by $\operatorname{Sym}_{d\times d}$ the set of real symmetric $d\times d$ matrices.
\item For a quadratic $2\times 2$ matrix $T$, we denote by $\overset{\circ}{T} = T-\frac{1}{2} \operatorname{tr}T \operatorname{Id}$ its traceless part.
\item For a function $f\in C^1(\R^2)$ we denote by $\nabla^\perp f = (\partial_2 f,-\partial_1 f)$ its orthogonal gradient.
\item For $d_1,d_2\in\N$ we write $f:\mathbb{T}^{d_1}\rightarrow\R^{d_2}$ for a function $f:\R^{d_1}\rightarrow\R^{d_2}$ defined on the full space that is periodic with period $1$ in all variables, i.e. $f(x+l e_k) = f(x)$ for all $k=1,\dots, d_1,$ $l\in\Z$.
\item For a periodic function $f$ as above, we denote $\int_{\mathbb{T}^{d_1}} f \,\dx = \int_{[0,1]^{d_1}} f\,\dx$, i.e. the integral over one periodic box.
\item $C_0^\infty(\tor;\R^d)=\lbrace f:\tor\rightarrow\R^d \text{ smooth}, \int_{\tor} f\,\dx = 0\rbrace$ is the space of smooth periodic functions on $\R^2$ with zero mean value on one periodic box.
\item For a function $g\in C^\infty(\tor;\R^d)$ and $\lambda\in\mathbb{N},$ we denote by $g_\lambda:\tor\rightarrow\mathbb{R}^d$ the $\frac{1}{\lambda}$ periodic function
$$g_\lambda (x) := g(\lambda x).$$
Notice that for every $l\in\mathbb{N}$, $s\in [1,\infty]$
$$\|D^l g_\lambda\|_{L^s(\tor)} = \lambda^l \|D^l g\|_{L^s(\tor)}.$$
\item $H^p(\R^2)$ is the real Hardy space, see Definition \ref{def: hardy space}.
\item $L^2_\sigma(\R^2) = \left\lbrace f\in L^2(\R^2): \dv f = 0 \text{ in distributions}\right\rbrace$ is the space of divergence-free vector fields in $L^2(\R^2)$.
\item For a function $f:[0,1]\times\R^2\rightarrow\R^d$ and $s\in[1,\infty]$, we write $\|\cdot\|_{C_tL^s_x}$ for the norm $\|f\|_{C_tL^s_x} = \max_{t\in[0,1]}\|f(t)\|_{L^s(\R^2)}$.
\item For any function $\varphi:\R\rightarrow\R$ with $\supp(\varphi)\subset (-\frac{1}{2},\frac{1}{2})$ and $\mu>1$ we write $\varphi_\mu$ for the periodic extension of the function $\mu^\frac{1}{2} \varphi(\mu \left(x-\frac{1}{2}\right))$, whose support is contained in intervalls of length $\frac{1}{\mu}$ centered around the points $\frac{1}{2}+\Z$. Note that
\begin{align*}
\|\varphi_\mu\|_{L^s(\mathbb{T})}=\mu^{\frac{1}{2} - \frac{1}{s}}\|\varphi\|_{L^s(\R)}, s\in [1,\infty],
\end{align*}
 and in particular $\|\varphi_\mu\|_{L^2(\mathbb{T})}=\|\varphi\|_{L^2(\R)}$.
\item Let $\lambda\in\N$, $f:\tor\rightarrow\R^d$. We will sometimes write $f_\lambda$ for the oscillating functions $f_\lambda(x) = f(\lambda x)$. On the other hand, for $f:\R\rightarrow\R$ with compact support, we will oftentimes write $f_\mu$ for its concentrated version. To avoid confusion, we will only use the parameter $\lambda$ for oscillations and $\mu$ (or $\mu_1,\mu_2$, respectively) for concentration.
\item We will use that for any concentrated and periodized function $\varphi_\mu(\lambda\cdot)\in C(\mathbb{T})$ it holds
\begin{align*}
\|D^l\varphi_\mu(\lambda\cdot)\|_{L^s(\mathbb{T})} = \lambda^l\mu^{l+\frac{1}{2}-\frac{1}{s}}.
\end{align*} 
\end{enumerate}

\section{Preliminaries}
We now provide the technicals tools that are needed for the proof of the Main Theorem \ref{thm:main}.

\subsection{An improved higher order inverse Laplacian}
As explained in the Introduction, we are working with perturbations of the form $f g_\lambda$, where $f\in C^\infty_c(\R^2)$ is compactly supported and $g\in C^\infty(\tor)$ is periodic, and we aim to write the perturbation $w_n$ introduced in \eqref{eq: introduction of wn}, that will be rigorously defined in Section \ref{sec: perturbations}, as $N$th order Laplacian. To this end, in Proposition \ref{prop: improved anti laplacian} below we construct a bilinear operator $\bl^N_M$ representing an inverse of $\Delta^N$ for functions of type $f g_\lambda$ up to a corrector of order $\lambda^{-M}$ (see Remark \ref{rem: inverse laplacian oscillations}) that in addition preserves the compact support of $fg_\lambda$. \\ For the next definition, we note that for $N\in\N$ and $f,g$ sufficiently smooth
\begin{align}\label{eq: higher order laplacian}
\Delta^N(fg) = \sum_{m_1+m_2+m_3=N}\sum_{|\alpha|=m_1}c_{(\alpha,m_2,m_3)}\partial^\alpha\Delta^{m_2} f  \partial^\alpha\Delta^{m_3}g
\end{align}
with some combinatorial constants $c_{(\alpha,m_2,m_3)}$.

\begin{dfn}
For $f,g\in C^\infty(\R^2)$,we define the commutator
\begin{align}\label{commutator}
\left[\Delta^N, f\right] g = \Delta^N(fg) - f\Delta^N g
\end{align}
and note that
\begin{align}\label{commutator sum}
\left[\Delta^N, f\right] g = \sum_{\substack{m_1+m_2+m_3=N\\m_3<N}}\sum_{|\alpha|=m_1}  c_{(\alpha,m_2,m_3)}\partial^\alpha\Delta^{m_2} f  \partial^\alpha\Delta^{m_3}g
\end{align}
with positive constants $c_{(\alpha,m_2,m_3)}$ from \eqref{eq: higher order laplacian}.
\end{dfn}
Before we can state and prove Proposition \ref{prop: improved anti laplacian}, we need to establish the following version of Calderón-Zygmund estimates:
\begin{lemma}\label{lemma: calderon zygmund}
Let $s\in [1,\infty]$, $m,l\in\N$ with $2l-1\leq m$. There exists a constant $C=C(s,m,l)>0$ such that for every $g\in C_0^\infty(\tor)$ it holds
\begin{align}\label{est: CZ endpoints}
\|g\|_{W^{m,s}(\tor)}\leq C \|\Delta^l g\|_{W^{m-(2l-1),s}(\tor)}.
\end{align}
In particular, 
\begin{align}\label{est: bound Delta-1}
\|\Delta^{-1} g\|_{W^{m+1,s}(\tor)} \leq \tilde{C}\|g\|_{W^{m,s}(\tor)}
\end{align}
for $s\in [1,\infty]$ and some $\tilde{C}(s,m)$.
\end{lemma}

\begin{proof}
For $s\in (1,\infty)$, this is a weaker version of the well-known Calderón-Zygmund estimate
\begin{align*}
\|g\|_{W^{2,p}(\tor)}\leq C \|\Delta g\|_{L^2(\tor)}
\end{align*}
applied $l$ times; see \cite{modena2020convex}, Lemma 3.2 for a proof of this estimate on the torus. For the endpoints $s\in \lbrace 1,\infty\rbrace$ one trades regularity against integrability; see the argument in \cite{modena2020convex}, Lemma 3.4. Estimate \eqref{est: bound Delta-1} follows immediately by picking $l=1$.
\end{proof}

\begin{prop}[Improved higher order inverse Laplacian]\label{prop: improved anti laplacian}
Let $N,M\in\N$. There exists bilinear operators 
\begin{align*}
\bl^N_M: C_c^\infty(\R^2)\times C_0^\infty(\tor)\rightarrow C_c^\infty(\R^2),\\
\vl^N_M: C_c^\infty(\R^2)\times C_0^\infty(\tor)\rightarrow C_c^\infty(\R^2)
\end{align*}
such that 
\begin{align*}
\Delta^N \bl^N_M(f,g) &= f g + \vl^N_M(f,g)
\end{align*}
with the following estimates for all $s\in [1,\infty]$
\begin{align}
\|\bl^N_M(f,g)\|_{W^{m,s}(\R^2)}&\leq C(N,M)  \|f\|_{C^{2NM+m}(\R^2)}\|\Delta^{-N}g\|_{W^{m,s}(\tor)},\label{est: bl1 no lambda}\\
\|\vl^N_M(f,g)\|_{W^{m,s}(\R^2)} &\leq C(N,M) \begin{cases}\|f\|_{C^{2N(M+1)+m}(\R^2)}\|\Delta^{-\frac{M}{2}}g\|_{W^{m,s}(\tor)}, \text{ if } M \text{ is even},\\
\|f\|_{C^{2N(M+1)+m}(\R^2)}\|\Delta^{-\frac{M+1}{2}}g\|_{W^{m+1,s}(\tor)}, \text{ if } M \text{ is odd}.
\end{cases}
\label{est: bl2 no lambda}
\end{align}
Furthermore, it holds
\begin{align}
\supp \bl^N_M(f,g), \supp \vl^N_M(f,g) \subseteq \supp f,\nonumber\\
\partial_t\bl^{N}_M(f,g) = \bl^{N}_M(\partial_t f, g) + \bl^{N}_M(f,\partial_t g),\nonumber\\
\partial_t\vl^{N}_M(f,g) = \vl^{N}_M(\partial_t f, g) + \vl^{N}_M(f,\partial_t g).\label{properties b v}
\end{align}
\end{prop}

\begin{rem}\label{rem: inverse laplacian oscillations}
We will use Proposition \ref{prop: improved anti laplacian} for fast oscillating functions $g_\lambda$: In that case, estimates \eqref{est: bl1 no lambda} and \eqref{est: bl2 no lambda} yield
\begin{align}
\|\bl^N_M(f,g_\lambda)\|_{W^{m,s}(\R^2)}&\leq C(N,M) \lambda^{-2N+m} \|f\|_{C^{2NM+m}(\R^2)}\|\Delta^{-N}g\|_{W^{m,s}(\tor)},\label{est: bl1}\\
\|\vl^N_M(f,g_\lambda)\|_{W^{m,s}(\R^2)} &\leq C(N,M) \lambda^{-M+m} \|f\|_{C^{2N(M+1)+m}(\R^2)}\|g\|_{W^{m,s}(\tor)}.\label{est: bl2}
\end{align}
Here, we used in \eqref{est: bl2} for even $M$ inequality \eqref{est: bound Delta-1}, while for odd $M$ we used first \eqref{est: CZ endpoints} with $l=1$ and then \eqref{est: bound Delta-1}:
\begin{align*}
\|\Delta^{-\frac{M+1}{2}}g\|_{W^{m+1,s}(\tor)}\leq C \|\Delta^{1-\frac{M+1}{2}} g\|_{W^{m,s}(\tor)}\leq C \|g\|_{W^{m,s}(\tor)}.
\end{align*}
\end{rem}

\begin{proof}[Proof of Proposition \ref{prop: improved anti laplacian}]
We recursively define bilinear operators $(b_k)_{k=0,\dots,M}$ that will cancel the bad parts of one another when their $N$th-order Laplacians $\Delta^N b_k$ are added up. We start with
\begin{align*}
b_0(f,g) &=f\Delta^{-N}g,\\
v_0(f,g) &= -fg
\end{align*}
which surely satisfy
\begin{align*}
\supp b_0(f,g), \supp v_0(f,g) &\subseteq \supp f ,\\
\partial_t b_0(f,g) &= b_0(\partial_t f, g) + b_0(f,\partial_t g),\\
\partial_t v_0(f,g) &= v_0(\partial_t f, g) + v_0(f,\partial_t g).
\end{align*}
The crucial obervation is that by \eqref{eq: higher order laplacian}
\begin{align}\label{laplace N}
\Delta^N b_0(f,g)= \Delta^N(f\Delta^{-N} g) &=\sum_{m_1+m_2+m_3=N}\sum_{|\alpha|=m_1} c_{(\alpha,m_2,m_3)} \partial^\alpha\Delta^{m_2} f  \partial^\alpha\Delta^{m_3-N}g\nonumber\\
&=\underbrace{f g}_{= -v_0(f,g)} +\underbrace{ \sum_{\substack{m_1+m_2+m_3=N\\m_3<N}}\sum_{|\alpha|=m_1} c_{(\alpha,m_2,m_3)}\partial^\alpha\Delta^{m_2} f \partial^\alpha\Delta^{m_3-N}g}_{=v_1(f,g), \text{ see below}}
\end{align}
with combinatorial constants $c_{(\alpha,m_2,m_3)}$ with $c_{(0,0,N)}=1$. The idea is to cancel the surplus terms of $\Delta^Nb_0(f,g)$,  by adding $b_1(f,g)$ which consists of the same terms as the surplus ones in  $\Delta^N b_0(f,g)$ but with  $\Delta^{-N}$ applied to the terms with $g$. Precisely, assuming $b_k$, $v_k$ are defined for some $k$, we define
\begin{align*}
b_{k+1}(f,g) &=  - \sum_{\substack{m_1+m_2+m_3=N\\m_3<N}}\sum_{|\alpha|=m_1} c_{(\alpha,m_2,m_3)} b_{k}(\partial^\alpha\Delta^{m_2} f,\partial^\alpha\Delta^{m_3-N}g),\\
v_{k+1}(f,g) &= - \sum_{\substack{m_1+m_2+m_3=N\\m_3<N}}\sum_{|\alpha|=m_1} c_{(\alpha,m_2,m_3)} v_k(\partial^\alpha\Delta^{m_2} f, \partial^\alpha\Delta^{m_3-N}g).
\end{align*}
In this way, the surplus terms of $\Delta^N b_{k}(f,g)$ are cancelled by those terms in $\Delta^N b_{k+1}(f,g)$ where $\Delta^N$ falls on the factor containing $g$. More precisely,
one can show by induction over $k$ and a direct computation that
\begin{align*}
\Delta^Nb_{k}(f,g) &= -v_{k}(f,g) + v_{k+1}(f,g)\text{ for all } k\in\N,
\end{align*}
from which it immediately follows
\begin{align*}
\Delta^N(\sum_{k=0}^{M}  b_{k}(f,g)) &= f g + v_{M+1}(f,g).
\end{align*}
Hence,  we define
\begin{align*}
\mathcal{B}^N_M(f,g) &=\sum_{k=0}^{M} b_{k}(f,g),\\
\vl^N_M(f,g) &= v_{M+1}(f,g),
\end{align*}
and therefore
\begin{align*}
\Delta^N\mathcal{B}^N_M(f,g) = fg + \vl^N_M(f,g).
\end{align*}
It is clear by induction that $\mathcal{B}^N_M(f,g)$ and $\vl^N_M(f,g)$ satisfy 
\begin{align}
\supp \mathcal{B}^N_M(f,g), \supp \vl^N_M(f,g) &\subseteq \supp f ,\nonumber\\
\partial_t \mathcal{B}^N_M(f,g) &= \mathcal{B}^N_M(\partial_t f,g) +\mathcal{B}^N_M(f,\partial_tg),\nonumber\\
\partial_t \vl^N_M(f,g) &= \vl^N_M(\partial_t f, g) + \vl^N_M(f,\partial_t g).\label{time derivative b and v}
\end{align}
For the estimate of $\mathcal{B}^N_M$, we show by induction on $k$ that
\begin{align}\label{est: gen product parts}
\|b_{k}(f,g)\|_{W^{m,s}(\R^2)} \leq C(N,k)\|f\|_{W^{2Nk+m,\infty}(\R^2)}\|\Delta^{-N}g\|_{W^{m,s}(\tor)}.
\end{align}\label{rem: number of derivatives}
That is true for $k=0$ by Hölder's inequality. Assuming the claim for some $k\in\N$, it follows by hypothesis
\begin{align*}
\|b_{k+1}(f,g)\|_{W^{m,s}(\R^2)}&\leq   \sum_{\substack{m_1+m_2+m_3=N\\m_3<N}}\sum_{|\alpha|=m_1} c_{(\alpha,m_2,m_3)} \|b_k\left(\partial^\alpha\Delta^{m_2} f, \partial^\alpha\Delta^{m_3-N}g\right)\|_{W^{m,s}(\R^2)}\\
 &\leq \sum_{\substack{m_1+m_2+m_3=N\\m_3<N}} \sum_{|\alpha|=m_1}c_{(\alpha,m_2,m_3)} \|\partial^\alpha\Delta^{m_2} f\|_{W^{2kN+m,\infty}(\R^2)} \\
&\hspace{5cm}\cdot\|\partial^\alpha\Delta^{m_3-N}(\Delta^{-N}g)\|_{W^{m,s}(\tor)}.
\end{align*}
It follows that the derivative in front of $f$ is of order $m_1+2m_2\leq 2N$ and therefore $$\|\partial^\alpha\Delta^{m_2} f\|_{W^{2Nk+m,\infty}(\R^2)}\leq \|f\|_{W^{2N(k+1)+m,\infty}(\R^2)}.$$ For the second factor, we use that $m_1+2m_3\leq 2N-1$ since $m_3<N$ and get 
\begin{align*}
 \|\partial^\alpha\Delta^{m_3-N}(\Delta^{-N}g)\|_{W^{m,s}(\tor)}\leq \|\Delta^{-N} (\Delta^{-N}g)\|_{W^{2N-1+m,s}(\tor)}\leq C(N)\|\Delta^{-N}g\|_{W^{m,s}(\tor)},
\end{align*} 
where we used Calderón-Zygmund estimates from Lemma \ref{lemma: calderon zygmund} in the second step, which proves \eqref{est: gen product parts} for $b_{k+1}$. We can now use \eqref{est: gen product parts} to  see that 
\begin{align}\label{est btilde}
\|\mathcal{B}^N_M(f,g)\|_{W^{m,s}(\R^2)} \leq C(N,M)\|f\|_{W^{2NM+m,\infty}(\R^2)}\|\Delta^{-N}g\|_{W^{m,s}(\tor)}.
\end{align}
In the same way, one can show by induction that for all $k\in\N$
\begin{align}
\|v_k(f,g)\|_{W^{m,s}(\R^2)}\leq C(N,k) \|f\|_{W^{2Nk+m,\infty}(\R^2)} \|\Delta^{-kN} g\|_{W^{(2N-1)k+m,s}(\tor)}. \nonumber
\end{align}
This yields by Calderón-Zygmund estimates from Lemma \ref{lemma: calderon zygmund} that for $l\leq \frac{(2N-1)k+1}{2}$
\begin{align*}
\|v_k(f,g)\|_{W^{m,s}(\R^2)}\leq C(N,k) \|f\|_{W^{2Nk+m,\infty}(\R^2)} \|\Delta^{l-kN} g\|_{W^{(2N-1)k-(2l-1)+m,s}(\tor)}.
\end{align*}
Choosing $l$ maximal, this yields
\begin{align*}
\|v_k(f,g)\|_{W^{m,s}(\R^2)} &\leq C(N,k) \begin{cases}\|f\|_{C^{2Nk+m}(\R^2)}\|\Delta^{-\frac{k-1}{2}}g\|_{W^{m,s}(\tor)}, \text{ if } k \text{ is odd},\\
\|f\|_{C^{2Nk+m}(\R^2)}\|\Delta^{-\frac{k}{2}}g\|_{W^{1+m,s}(\tor)}, \text{ if } k \text{ is even}.
\end{cases}
\end{align*}
This shows also the estimate for $\vl^N_M= v_{M+1}(f,g)$. 
\end{proof}

\subsection{Hardy spaces}\label{subsec: hardy spaces}
In this subsection, we first collect some basic properties of real Hardy spaces $H^p(\R^d)$ and then prove an estimate for finite collections of functions whose sum is in $H^p(\R^d)$ but the functions themselves are not necessarily.
\subsubsection{Basics on Hardy spaces}
\begin{dfn}[Hardy spaces on $\R^d$]\label{def: hardy space}
 Let $\Psi\in \mathcal{S}(\R^d)$ be a Schwartz function with $\int_{\R^2}\Psi(x)\,\dx \neq 0$ and let $\Psi_\zeta(x) = \frac{1}{\zeta^d}\Psi(\frac{x}{\zeta})$. For any $f\in \mathcal{S}'(\R^d),$ we define the  radial maximal function
\begin{align}
\label{eq_def_hardy_norm}
m_\Psi f(x) = \sup_{\zeta>0}|f\ast\Psi_\zeta(x)|.
\end{align}
Let $0<p<\infty$. The real Hardy space $H^p(\R^d)$ is defined as the space of tempered distributions
\begin{align*}
H^p(\R^d) = \left\lbrace f\in \mathcal{S}'(\R^d) : m_\Psi f \in L^p(\R^d)\right\rbrace
\end{align*}
and we write
\begin{align*}
\|f\|_{H^p(\R^d)} = \|m_\Psi f\|_{L^p(\R^d)}.
\end{align*}
Note that $\|\cdot\|_{H^p(\R^d)}$ is only a quasinorm, in particular the triangle inequality does not hold, instead $\|f-g\|_{H^p(\R^d)}\leq K(\|f\|_{H^p(\R^d)} + \|g\|_{H^p(\R^d)})$ for some $K>1$. 
The definition of $H^p(\R^d)$ does not depend on the choice of the function $\Psi$ and the quasinorms are equivalent. For $p>1$, the space $H^p(\R^d)$ coincides with the Lebesgue space $L^p(\R^d)$. For $p\leq 1$, $H^p(\R^d)$ is a complete metric space with the metric given by $d(f,g)=\|f-g\|^p_{H^p(\R^d)}$ and the inclusion $H^p(\R^d)\hookrightarrow \mathcal{S}'(\R^d)$ is continuous, see \cite{grafakos2009modern}, Proposition 6.4.10.
\end{dfn}

\begin{dfn}[Hardy space atoms]\label{dfn: hardy atoms}
For $p\leq 1$, a Hardy space atom is a measurable function $a$ with the following properties:
\begin{enumerate}[(i)]
\item $\supp a\subset B$ for some ball $B$,
\item $|a|\leq |B|^{-\frac{1}{p}}$
\item $\int_B x^\beta a(x)\, \dx = 0$ for all multiindices $\beta$ with $|\beta|\leq d(p^{-1}-1)$.
\end{enumerate}
\end{dfn}

\begin{lemma}[Estimate for Hardy space atoms]\label{lemma: hardy atoms}
There is a uniform constant $C$ such that for all atoms $a$ it holds
\begin{align*}
\|a\|_{H^p(\R^d)}\leq C.
\end{align*}
\end{lemma}
\begin{proof}
We refer to \cite{stein1993harmonic}, see 2.2 in Chapter III.2.
\end{proof}

\begin{rem}\label{rem: hardy atoms}
We will use that for a function $f$ satisfying (iii) in Definition \ref{dfn: hardy atoms} with support in a ball $B$, we have $f\in H^p(\R^d)$ and by \mbox{Lemma \ref{lemma: hardy atoms}}
\begin{align}\label{est scaled atom}
\|f\|_{H^p(\R^d)} \leq C |B|^\frac{1}{p}\|f\|_{L^\infty(\R^d)}.
\end{align}

\end{rem}

\subsubsection{Estimate for sums in Hardy spaces}
As pointed out in the Introduction, we will define perturbations consisting of different parts in such a way that their vorticities are in $H^p(\R^2)$, but the vorticities of their individual parts are not necessarily. In order to deal with such perturbations, we prove the following proposition. 

\begin{prop}[Hardy estimate for sums]\label{Hardy for sum}
Let $J\in\N$, $(y_j)_{1\leq j\leq J}\subset\R^d$ and $(f_j)_{1\leq j\leq J}\subset L^\infty(\R^d)$ be compactly supported functions with $\supp f_j\subset B_{\varepsilon_j}(y_j)$, $\varepsilon_j>0$, such that $\sum_{j=1}^J f_j\in H^p(\R^d)$. Furthermore, let $N\in\N$ be maximal such that $N\leq d(\frac{1}{p}-1)$ and let $R>0$ such that \mbox{$\bigcup_j B_{\varepsilon_j}(y_j)\subset B_R(y)$} for some $y\in\R^d$. Then
\begin{align}\label{est hardy sum}
\|\sum_{j=1}^J f_j\|_{H^p(\R^d)}^p &\leq  C(p) \sum_{j=1}^J\varepsilon_j^d  \|f\|_{L^\infty(\R^d)}^p + \varepsilon_j^{d-p(d+N)} \max\left\lbrace R^{Np}, 1\right\rbrace\max_{|\alpha|\leq N}\left|\int x^\alpha f_j(x) \,\dx\right|^p\nonumber\\
&\hspace{0,3cm} +C(p) \max\lbrace R^{d-pd}, R^{d-p(d+N)}\rbrace \sum_{j=1}^J \left(1+\log\left(\frac{R}{\varepsilon_j}\right)\right)\max_{|\alpha|\leq N}\left|\int x^\alpha f_j(x) \,\dx\right|^p.
\end{align}
In case $N<d(\frac{1}{p}-1)$, i.e. if $\frac{d}{p}\notin\N$, \eqref{est hardy sum} holds without the term $\log\left(\frac{R}{\varepsilon_j}\right)$ in the second line.
\end{prop}

\begin{cor}If $R=1$ and $\varepsilon_j< \operatorname{e}^{-1}$ we can use that  $d-p(d+N)\geq 0$, i.e. \mbox{$\varepsilon_j^{d-p(d+N)}\leq 1\leq |\log(\varepsilon_j)|$}, therefore estimate \eqref{est hardy sum} becomes
\begin{align}\label{est hardy sum small eps}
\|\sum_{j=1}^J f_j\|_{H^p(\R^d)}^p&\leq C(p) \sum_{j=1}^J\varepsilon_j^d \|f_j\|^p_{L^\infty(\R^d)} + |\log(\varepsilon_j)| \max_{|\alpha|\leq N}\left|\int x^\alpha f_j(x) \,\dx\right|^p.
\end{align}
\end{cor}
To prove Proposition \ref{Hardy for sum}, we need some auxiliary functions  given by Lemma \ref{lemma: improved anti laplacian} below. We quickly introduce some\\ 
\noindent \textbf{Notation.}  It is possible to define a total order $\preceq$ on the set $\N^d$ of multi-indexes in such a way that, if 
\begin{equation*}
|\alpha| < |\beta| \quad \Longrightarrow \quad \alpha \preceq \beta. 
\end{equation*}
\noindent In particular, if $\alpha \leq \beta$, then $\alpha \preceq \beta$. Using such order, any family of real numbers $\{M_\alpha\}_{|\alpha| \leq N}$ can be written as a vector in $\R^{\binom{N+d}{N}}$, $\binom{N+d}{N}$ being the cardinality of the set $\{ \alpha \in \N^d \, : \, |\alpha| \leq N\}$ of all multi-indexes up to order $N$. Similarly, any family  $\{M_{\alpha \beta}\}_{|\alpha|, |\beta| \leq N}$ of real numbers indexed by pairs of multi-indexes of order up to $N$,  can be written as a square matrix in $\R^{\binom{N+d}{N} \times \binom{N+d}{N}}$. For example, for $d = 2$, one possible choice for the ordering $\preceq$  is
\begin{equation*}
(0,0), (0,1), (1,0), (0,2), (1,1), (2,0), (0,3), (1,2), (2,1), (3,0), \dots
\end{equation*}
W.r.t. this order a family of real numbers $\{M_{\alpha \beta}\}_{|\alpha|, |\beta| \leq N}$ (for $N =1$) can be written in matrix form as
\begin{equation*}
\left(
\begin{matrix}
M_{(0,0), (0,0)} & M_{(0,0), (0,1)} & M_{(0,0), (1,0)} \\
M_{(0,1), (0,0)} & M_{(0,1), (0,1)} & M_{(0,1), (1,0)} \\
M_{(1,0), (0,0)} & M_{(1,0), (0,1)} & M_{(1,0), (1,0)} \\
\end{matrix}
\right).
\end{equation*}
The following lemma shows the existence of a linear operator that assigns to each given tuple \mbox{$m\in \R^{\binom{N+d}{N}}$} of momenta up to order $N$ and to each given point $y\in\R^d$ in space a smooth function with that exact prescribed momenta and support in an arbitrarily small ball around $y$.
\begin{lemma}\label{lemma: improved anti laplacian}
For every $N \in \N$ there is a degree $N$ polynomial $p_N$ on $\R$ such that the following holds. For all $\e \in (0,1]$ and $y \in \R^d$, there is a linear operator
\begin{equation*}
\Psi_{\e, y} : \R^{\binom{N+d}{N}} \to C^\infty_c(B_\e(y)), \qquad m \mapsto \Psi_{\e, y, m}
\end{equation*}
such that 
\begin{equation}
\label{eq:momenta}
\int x^\alpha \Psi_{\e, y, m}(x) dx = m_\alpha \text{ for all } |\alpha| \leq N
\end{equation}
and
\begin{equation}
\label{eq:momenta-bound}
\|\Psi_{\e, y,m}\|_{L^\infty} \leq \frac{p_N(|y|)}{\e^{d + N}} |m|.
\end{equation}
\end{lemma}
\begin{proof}
Let $\varphi \in C^\infty_c(B_1(0))$ with $\int \varphi = 1$. The notation $C_N$ denotes a constant which depends on $N$, the fixed function $\varphi$, but not on $y$, $\e$ and $m$. Now given $\e, y, m$, we define
\begin{equation*}\
 \Psi_{\e, y, m}(x) = \sum_{|\zeta| \leq N} a_\zeta \left( \partial^\zeta \varphi \right) \left( \frac{x - y}{\e} \right)
\end{equation*}
for coefficients $a = (a_\zeta)_{|\zeta| \leq N}$ to be suitably chosen. Observe that $\Psi_{\e, y, m}(x) \in C^\infty_c(B_\e(y))$ and $\|\Psi_{\e, y, m}(x)\|_{L^\infty} \leq C_N |a|$. We have thus to show that we can choose the vector  $a$ of coefficients such that \eqref{eq:momenta} and \eqref{eq:momenta-bound} are satisfied and $a$ depends linearly on $m$. We compute
\begin{equation*}
\begin{aligned}
\int x^\alpha \Psi_{\e, y, m}(x) \dx 
& = \sum_{|\zeta| \leq N} a_\zeta \int x^\alpha  \left( \partial^\zeta \varphi \right) \left( \frac{x - y}{\e} \right) \dx \\
\text{(changing variable)}
& = \sum_{|\zeta| \leq N} a_\zeta \int (y+\e z)^\alpha  \partial^\zeta \varphi (z) \e^d \dz \\
& = \sum_{|\zeta| \leq N} \sum_{\beta \leq \alpha} a_\zeta \binom{\alpha}{\beta} y^{\alpha - \beta} \e^{|\beta| + d} \int z^\beta \partial^\zeta \varphi (z) \dz \\
\text{(writing $\sum_{|\gamma| \leq N} \delta_{\beta \gamma} \int z^\gamma \partial^\zeta \varphi(z) \dz$ )}
& = \sum_{|\zeta| \leq N}  \sum_{\beta \leq \alpha} \sum_{|\gamma| \leq N} a_\zeta \underbrace{\binom{\alpha}{\beta} y^{\alpha - \beta}}_{=:A_{\alpha \beta}} \underbrace{\e^{|\beta| + d} \delta_{\beta \gamma}}_{=: B_{\beta \gamma}} \underbrace{\int z^\gamma \partial^\zeta  \varphi (z) \dz}_{=: M_{\gamma \zeta}} \\
& = \sum_{|\zeta| \leq N}  \sum_{|\beta| \leq N} \sum_{|\gamma| \leq N} A_{\alpha \beta} B_{\beta \gamma} M_{\gamma \zeta} a_\zeta,
\end{aligned}
\end{equation*} 
where we have defined
\begin{equation*}
A_{\alpha \beta} :=
\begin{cases}
\binom{\alpha}{\beta} y^{\alpha - \beta}, & \text{if } \beta \leq \alpha, \\
0, & \text{otherwise},
\end{cases}
\qquad
B_{\beta \gamma} := \e^{|\beta| + d} \delta_{\beta \gamma}, 
\qquad
M_{\gamma \zeta} := \int z^\gamma \partial^\zeta \varphi(z) \dz.
\end{equation*}
Notice that in the last inequality above the summation over $\beta$ runs for all $|\beta| \leq N$ since $A_{\alpha \beta} = 0$ if $\beta\leq \alpha$ is  not satisfied. In matrix notation we thus have
\begin{equation*}
\left( \int x^\alpha \Psi_{\e, y, m}(x) \dx \right)_{|\alpha| \leq N} = A B M a. 
\end{equation*}
Hence \eqref{eq:momenta} is satisfied if we can define $a$ such that $m = A B M a$. Observe that 
\begin{itemize}
\item $A$ is lower triangular (by definition) and its diagonal elements are all equal to $1$. Hence $A$ is invertible. Standard combinatorial computation in linear algebra shows also that
\begin{equation*}
A^{-1}_{\alpha \beta} =
\begin{cases}
(-1)^{|\alpha - \beta|}\binom{\alpha}{\beta} y^{\alpha - \beta}, & \text{if } \beta \leq \alpha, \\
0, & \text{otherwise},
\end{cases}
\end{equation*}
and thus, in particular,
\begin{equation*}
|A^{-1}| \leq C_N (1 + |y|^N) =: p_N(|y|).
\end{equation*}

\item $B$ is diagonal, its inverse is given by $B^{-1}_{\beta \gamma} = \e^{- (|\beta| + d)} \delta_{\beta \gamma}$ and thus, since $\varepsilon\in (0,1]$
\begin{equation*}
|B^{-1}| \leq \e^{- (N+d)}.
\end{equation*}
\item $M$ is lower triangular, because, by integration by parts, $\int z^\gamma \partial^\zeta \varphi(z) \,\dz\neq 0$ only if $\gamma \geq \beta$. Moreover, again by integration by parts and using that $\int \varphi = 1$, the diagonal elements are all nonzero. Hence $M$ is invertible and its inverse is bounded by some constant $|M^{-1}| \leq C_N$. 

\end{itemize}
Therefore defining $a:= M^{-1} B^{-1} A^{-1} m$, we obtain \eqref{eq:momenta} and the linear dependence on $a$ (and thus of $\Psi_{\e, y, m}$) on $m$. Moreover (possibly increasing the value of the constant $C_N$),
\begin{equation*}
|a| \leq |M^{-1}| | B^{-1}| | A^{-1}| |m| \leq \frac{p_N(|y|)}{\e^{N + d}} |m|,
\end{equation*}
thus concluding the proof of the proposition. 
\end{proof}

\begin{proof}[Proof of Proposition \ref{Hardy for sum}]
We show the estimate for $R=1$ and then argue by scaling. For each $j$, we choose a sequence of balls with shrinking radii
\begin{align*}
1=\varepsilon_j^0 &> \varepsilon_j^1 >\dots > \varepsilon_j^{K_j} = \varepsilon_j
\end{align*}
where $K_j\in\N$ is minimal with $2^{-K_j}\leq \varepsilon_j$ and we set
\begin{align*}
\varepsilon_j^n&= 2^{-n} \text{ for } n< K_j,\\
\varepsilon_j^{K_j} &= \varepsilon_j
\end{align*}
and choose points $y=y_j^0, y_j^1,\dots, y_j^{K_j}=y_j$ such that
\begin{align*}
B_1(y)= B_{\varepsilon_j^0}(y_j^0)&\supseteq   B_{\varepsilon_j^1}(y_j^1)\supseteq \dots \supseteq B_{\varepsilon_j^{K_j}}(y_j^{K_j})= B_{\varepsilon_j}(y_j)
\end{align*}
which is possible since $B_{\varepsilon_j}(y_j)\subset B_1(y)$ by assumption. We cannot directly use estimate \eqref{est scaled atom} on the single summands and use the triangle inequality for $\|\cdot\|_{H^p(\R^d)}^p$ because we do not assume that $f_j\in H^p(\R^d)$. Therefore,  we decompose $\sum_{j=1}^J f_j$ into functions in $H^p(\R^d)$ using a telescopic sum: Let $m^j\in\R^{\binom{N+d}{N}}$ given by $m^j_\alpha = \int x^\alpha f_j \,\dx$. Using  the functions from Lemma \ref{lemma: improved anti laplacian},  we decompose
\begin{align*}
\sum_{j=1}^J f_j &= \sum_{j=1}^J f_j - \Psi_{\varepsilon_j, y_j, m^j} +  \Psi_{\varepsilon_j, y_j, m^j} \\
&=  \sum_{j=1}^J f_j - \Psi_{\varepsilon_j, y_j, m^j} \\
&\hspace{0,3cm} + \sum_{j=1}^J\sum_{k=0}^{K_j-1} \Psi_{\varepsilon_j^{k+1}, y_j^{k+1}, m^j} -  \Psi_{\varepsilon_j^{k}, y_j^{k}, m^j}\\
& \hspace{0,3cm}+ \sum_{j=1}^J \Psi_{1,y,m^j}. 
\end{align*}
For the last line, we notice that by linearity of $\Psi_{1, y}$ in $m$ 
\begin{align*}
 \sum_{j=1}^J \Psi_{1,y,m^j} &= \Psi_{1,y,\sum_j m^j}= 0
\end{align*}
because $(\sum_j m^j)_\alpha = \int x^\alpha \left(\sum_{j=1}^J f_j\right)\,\dx=0$ because of  $\sum_{j=1}^J f_j\in H^p(\R^d)$, which means that all momenta up to order $N$ vanish. Also, from the properties of $\Psi_{\varepsilon_j, y_j, m^j}$, it is clear that \mbox{$f_j - \Psi_{\varepsilon_j, y_j, m^j}\in H^p(\R^d)$}  for every $1\leq j\leq J$ since this sum is supported in $B_{\varepsilon_j}(y_j)$ and all momenta up to order $N$ vanish by Lemma \ref{lemma: improved anti laplacian}. With the same reasoning, we see that $\Psi_{\varepsilon_j^{k+1}, y_j^{k+1}, m^j} -  \Psi_{\varepsilon_j^{k}, y_j^{k}, m^j} \in H^p(\R^d)$ for every $1\leq j\leq J$ and $0\leq k\leq K_j-1$. Therefore, we can estimate
\begin{align}
\|\sum_{j=1}^J f_j\|_{H^p(\R^d)}^p &\leq  \sum_{j=1}^J\|f_j - \Psi_{\varepsilon_j,y_j,m^j}\|_{H^p(\R^d)}^p\nonumber\\
&\hspace{0,3cm} + \sum_{j=1}^J\sum_{k=0}^{K_j-1}\|\Psi_{\varepsilon_j^{k+1}, y_j^{k+1}, m^j} -  \Psi_{\varepsilon_j^{k}, y_j^{k}, m^j}\|_{H^p(\R^d)}^p.\label{est: telescopic}
\end{align}
We estimate the first summand in the previous inequality by \eqref{est scaled atom} and  by the estimate in Lemma \ref{lemma: improved anti laplacian}: We have
\begin{align*}
\|f_j - \Psi_{\varepsilon_j,y_j,m^j}\|^p_{H^p(\R^d)} &\leq  C\varepsilon_j^d \left(\|f_j\|_{L^\infty(\R^2)} + \frac{1}{\varepsilon^{d+N}}\max_{|\alpha|\leq N}\int x^\alpha f_j(x)\,\dx\right)^p.
\end{align*}
Considering the sums in the second line of \eqref{est: telescopic}, we have for each $1\leq j\leq J$ and $0\leq k\leq K_j-1$ by the estimate in Lemma \ref{lemma: improved anti laplacian}
\begin{align*}
\|\Psi_{\varepsilon_j^{k+1}, y_j^{k+1}, m^j} -  \Psi_{\varepsilon_j^{k}, y_j^{k}, m^j}\|_{L^\infty(\R^d)} \leq C\frac{1}{(\varepsilon_j^{k+1})^{d+N}}  \max_{|\alpha|\leq N} \left|\int x^\alpha f_j(x)\,\dx\right|.
\end{align*}
Therefore, since $\Psi_{\varepsilon_j^{k+1}, y_j^{k+1}, m^j} -  \Psi_{\varepsilon_j^{k}, y_j^{k}, m^j}$ is supported in $B_{\varepsilon_j^k}(y_j^k)$, by the estimate for multiples of atoms 
\begin{align*}
&\sum_{k=0}^{K_j-1}\|\Psi_{\varepsilon_j^{k+1}, y_j^{k+1}, m^j} -  \Psi_{\varepsilon_j^{k}, y_j^{k}, m^j}\|^p_{H^p(\R^d)}\leq  C\sum_{k=0}^{K_j-1} \frac{(\varepsilon_j^{k})^d }{(\varepsilon_j^{k+1})^{p(d+N)}}\max_{|\alpha|\leq N} \left|\int x^\alpha f_j(x)\,\dx\right|^p.
\end{align*}
We use  that $\varepsilon_j^{K_j}\geq 2^{-K_j}$, and the minimality of $K_j$ also yields 
\begin{align*}
K_j\leq C\left(1+\log\left(\frac{1}{\varepsilon_j}\right)\right).
\end{align*} 
Exploiting also the fact that
\begin{align*}
\begin{cases}
N&< d(\frac{1}{p}-1), \frac{d}{p}\notin\N,\\
N&= d(\frac{1}{p}-1), \frac{d}{p}\in \N,
\end{cases}
\end{align*}
and the geometric sum formula, we have
\begin{align}
\sum_{k=0}^{K_j-1}  \frac{(\varepsilon_j^{k})^d }{(\varepsilon_j^{k+1})^{p(d+N)}} &\leq \sum_{k=0}^{K_j-1} \frac{2^{-kd}}{2^{-(k+1)p(d+N)}}\leq C(p)  \sum_{k=0}^{K_j-1} 2^{-k(d-p(d+N))}\nonumber\\
&\leq C(p)  \begin{cases} 1, &\frac{d}{p}\notin \N,\\
K_j, &\frac{d}{p}\in\N 
\end{cases} \leq   C (p) \begin{cases} 1, &\frac{d}{p}\notin \N\\
 1+\log\left(\frac{1}{\varepsilon_j}\right), &\frac{d}{p}\in\N .
\end{cases}\label{est case distinction log}
\end{align}
This yields the claim for $R=1$. Estimate \eqref{est case distinction log} shows that we can omit the term including the logarithm if $\frac{d}{p}\notin \N$. For $R\neq 1$, we consider the functions $\tilde{f}_j(x) := f_j(Rx)$, which are supported on balls with radii $\tilde{\varepsilon}_j = \frac{\varepsilon_j}{R}<1$ around $\frac{y_j}{R}$ and $\supp \cup_j \tilde{f}_j\subseteq B_1(\frac{y}{R})$ and use the estimate that we have already established, noting that $\|\sum_{j=1}^J \tilde{f}_j\|_{H^p(\R^2)}^p = R^{-d}\|\sum_{j=1}^J f_j\|_{H^p(\R^2)}^p$ and
\begin{align*}
\left|\int x^\alpha \tilde{f}_j(x) \,\dx \right| \leq \begin{cases} \frac{1}{R^d} \left|\int x^\alpha f_j(x) \,\dx\right| &\text{ if } R>1,\\
\frac{1}{R^{d+N}} \left|\int x^\alpha f_j(x) \,\dx\right| &\text{ if } R<1.
\end{cases}
\end{align*}
\end{proof}

\subsection{Antidivergence operators}
In this section, we present various antidivergence operators. As pointed out in the introduction, one main obstruction to performing convex integration techniques on the full space $\R^d$ compared to the flat torus is the absence of a bounded operator $\dv^{-1}:L^1(\R^d;\R^d)\rightarrow L^1(\R^d;\operatorname{Sym}_{d\times d})$. We present in Proposition \ref{thm: isett} a result by Isett and Oh \cite{isett2016fullspace} that yields a bounded antidivergence operator on the set of vector fields with vanishing linear and angular momentum and use that operator to define an improved bilinear antidivergence operator in Lemma \ref{antidivergence}.
\begin{lemma}[Standard antidivergence on the torus]\label{standard antidivergence}
There exists a linear operator $$\antidvp: C^\infty_0(\tor;\R^2)\rightarrow C^\infty_0(\tor;\operatorname{Sym}_{2\times 2})$$ such that $\dv \antidvp u =u$ and
\begin{align*}
\|\nabla^l\antidvp u\|_{L^s(\tor)}&\leq C(s)\|\nabla^lu\|_{L^s(\tor)} ,\\
\|\nabla^l\antidvp u_\lambda\|_{L^s(\tor)}&\leq \frac{C(s)}{\lambda^{1-l}}\|\nabla^l u\|_{L^s(\tor)}\text{ for all } l,\lambda\in\N,s\in [1,\infty].
\end{align*}
\end{lemma}
For the proof see Proposition 4 in \cite{burczak-mod-sze21}.

\begin{prop}[Isett-Oh]\label{thm: isett}
There exists a constant $C>0$ such that the following holds. Let $X$ be the subspace of $C_c^\infty(B_1(0);\R^d)$ consisting of those vector fields with vanishing linear and angular momentum
\begin{align*}
X := \left\lbrace f\in C_c^\infty(B_1(0);\R^d): \int f\, \dx =0, \int x_j f_l - x_l f_j \,\dx  = 0 \text{ for all } i,j=1,\dots,d\right\rbrace.
\end{align*} 
There exists a linear antidivergence operator $\antidvx: X\rightarrow C_c^\infty(B_1(0); \operatorname{Sym}_{d\times d})$  with $$\supp \antidvx f \subseteq B_1(0),$$ for all $f\in X$, that solves
\begin{align}\label{divergence}
\dv\antidvx f  = f
\end{align}
and satisfies for all $s\in [1,\infty]$
\begin{align}\label{est: Risett}
\|\antidvx f\|_{L^s(B_1(0))} \leq C \|f\|_{L^s(B_1(0))}.
\end{align}
Furthermore, if $f = f(t,x)$ depends on time, then $\antidvx f(t,x) := \antidvx(f(t,\cdot))(x)$ is also smooth in time. 
\end{prop}
For a detailed proof of Proposition \ref{thm: isett}, see Section 11 in \cite{isett2016fullspace}.

\begin{dfn}\label{proj a}
Let $\iota\in C^\infty_c(B_1(0); \R), \eta\in C^\infty_c(B_1(0); \R^2)$ with
\begin{align*}
\int \iota \,\dx &= 1, \int x_1 \iota\,\dx = \int x_2 \iota\,\dx = 0,\\
\int \eta\,\dx &= 0, \int x_1\eta_2 - x_2\eta_1 \,\dx = 1.
\end{align*}
With $\iota$ and $\eta$ from Definition \ref{proj a} we define for $f\in C_c^\infty(B_1(0);\R^2)$ a projection onto $X$ as
\begin{align*}
\mathbb{P}_X f = f- \iota\int f\,\dx - \eta \int\left( x_1f_2(x) - x_2f_1(x) \right)\,\dx.
\end{align*}
Since $f$ is compactly supported, it holds for all $s\in [1,\infty]$
\begin{align}\label{est: px}
\|\px f\|_{L^s(\R^2)}\leq C(s) \|f\|_{L^s(\R^2)}.
\end{align}
\end{dfn}

With the standard antidivergence operator on $\tor$ and Isett and Oh's operator, we will define an improved antidivergence operator for functions of the form $fu_\lambda$, $f\in C^\infty_c(\R^2)$, $u\in C^\infty_0(\tor;\R^2)$, on the full space in the next lemma.
\begin{lemma}[Improved antidivergence operators]\label{antidivergence}
\begin{enumerate}[(i)]
\item There exists a bilinear operator $$\Rc: C_c^\infty(B_1(0);\R)\times C_0^\infty(\tor;\R^2)\rightarrow C_c^\infty(B_1(0);\operatorname{Sym}_{2\times 2})$$ such that  it holds
\begin{align*}
\dv \Rc(f,u) = \px(fu)
\end{align*}
with
\begin{align*}
\|\Rc(f,u)\|_{L^1(\R^2)} &\leq C  \|f\|_{W^{1,\infty}(\R^2)} \|\antidvp u\|_{L^1(\tor)}.
\end{align*}
\item There exists a bilinear operator $$\tilde{\Rc}: C_c^\infty(B_1(0);\R^2)\times C_0^\infty(\tor;\operatorname{Sym}_{2\times 2})\rightarrow C_c^\infty(B_1(0);\operatorname{Sym}_{2\times 2})$$ such that it holds
\begin{align*}
 \dv \tilde{\Rc}(f,T) = \px(Tf)
\end{align*}
with
\begin{align*}
\|\tilde{\Rc}(f,T)\|_{L^1(\R^2)} &\leq  C \| f\|_{W^{1,\infty}(\R^2)} \|\antidvp T\|_{L^1(\tor)}
\end{align*}
where, by a slight abuse of notation, we define
\begin{align*}
\dv^{-1}_{\mathbb{T}^2} T= \sum_{k=1,2} \dv^{-1}_{\mathbb{T}^2}(Te_k).
\end{align*}
\end{enumerate}
\end{lemma}
\begin{proof}
Let us  define
\begin{align*}
\Rc:C_c^\infty(B_1(0);\R)\times C_0^\infty(\tor;\R^2)&\rightarrow C_c^\infty(B_1(0);\operatorname{Sym}_{2\times 2}),\\
\Rc(f,u) &= f\antidvp u -\antidvx\left(\mathbb{P}_X( (\antidvp u) \nabla f)\right).
\end{align*}
It is not difficult to check that $\dv \Rc(f,u) = \mathbb{P}_X(fu),$ since with $\iota$ and $\eta$ as in Definition  \ref{proj a}, we have
\begin{align*}
\dv \Rc(f,u) &=   \dv\left(f\antidvp u - \antidvx\left(\mathbb{P}_X\left((\antidvp u) \nabla f\right)\right)\right)\\
&= fu + (\antidvp u)\nabla f  -\px(  (\antidvp u) \nabla f)\\
&= fu + \iota\int (\antidvp u)\nabla f \,\dx \\
&\hspace{0,3cm}+ \eta \int x_1 \left((\antidvp u)\nabla f\right)_2 - x_2 \left((\antidvp u)\nabla f\right)_1\,\dx\\
&= \px(fu)
\end{align*}
through integration by parts in the last equality.
For the second operator, we simply set for $f\in C_c^\infty(B_1(0);\R^2)$, $T\in C_0^\infty(\tor;\operatorname{Sym}_{2\times 2})$
\begin{align*}
\tilde{\Rc}(f,T) = \sum_{k=1,2} \Rc(f_k, Te_k).
\end{align*}
The estimates follow  from the ones for $\antidvp$ from Lemma \ref{standard antidivergence}, $\antidvx$ from Proposition \ref{thm: isett} and the estimate \eqref{est: px} for $\px$.
\end{proof}

\begin{rem}\label{rem: improved antidiv}
In particular, for fast oscillating $u_\lambda\in C^\infty_0(\tor)$ and $f\in C^\infty_c(B_1(0))$, we see that $\supp\Rc(f,u_\lambda)\subseteq B_1(0)$, i.e. we keep the compact support of $fu_\lambda$. Using that $\antidvp$ is bounded in $L^1(\tor)$,
\begin{align*}
\|\Rc(f,u_\lambda)\|_{L^1(\R^2)} &\leq \frac{C}{\lambda} \|f\|_{W^{1,\infty}(\R^2)} \|u\|_{L^1(\tor)}
\end{align*}
and the same holds for $\tilde{\Rc}$.
\end{rem}

\subsection{Further helpful tools}
We close this preliminary section with a useful estimate for functions of the form  $fg_\lambda$, where \mbox{$f\in C_c^\infty(\R^2),$}  $g\in C^\infty(\tor)$ and some identities for products of "orthogonally oriented" functions.
\begin{lemma}[Improved Hölder]\label{lemma: improved hoelder}
Let  $\lambda\in\mathbb{N},$ $f:\R^2\rightarrow\R$ smooth with $\supp f\subset B_1(0)$ and $g:\mathbb{T}^2\rightarrow\R$ smooth. Then it holds for all $s\in [1,\infty]$
\begin{align*}
\|fg_\lambda\|_{L^s(\R^2)}\leq \|f\|_{L^s(\R^2)}\|g\|_{L^s(\mathbb{T}^2)} + \frac{C(s)}{\lambda^\frac{1}{s}}\|f\|_{C^1(\R^2)}\|g\|_{L^s(\mathbb{T}^2)}.
\end{align*}
\end{lemma}

\begin{proof}
We refer to Lemma 2.1 in \cite{modena2018non}. For this estimate it is crucial that $f$ is compactly supported. Note also that the size of $\supp f$ enters the estimate.
\end{proof}



\begin{lemma}\label{lemma: div of matrix}
Let $f,g\in C^1(\R)$. For any vector $\xi\neq 0\in \R^2$ it holds
\begin{align*}
\dv\left( f(\xi\cdot x) g(\xi^\perp\cdot x) \frac{\xi}{|\xi|}\otimes  \frac{\xi}{|\xi|}\right) &= f'(\xi\cdot x) g(\xi^\perp\cdot x) \xi,\\
\dv\left(f(\xi\cdot x) g(\xi^\perp\cdot x) \frac{\xi^\perp}{|\xi|}\otimes \frac{\xi}{|\xi|} \right) &= f'(\xi\cdot x) g(\xi^\perp\cdot x)\xi^\perp,\\
\dv\left(f(\xi\cdot x) g(\xi^\perp\cdot x)\xi\right) &
= f'(\xi\cdot x) g(\xi^\perp\cdot x)|\xi|^2,\\
\Delta\left(f(\xi\cdot x) g(\xi^\perp\cdot x)\right) &= f''(\xi\cdot x) g(\xi^\perp\cdot x)|\xi|^2 + f(\xi\cdot x) g''(\xi^\perp\cdot x)|\xi|^2.
\end{align*}
\end{lemma}
\begin{proof}
The proof is trivial.
\end{proof}

\section{Main Proposition}
In this section we present the main proposition that is the key to prove Theorem \ref{thm:main}. To this end, we first introduce the Reynolds-defect-equation:
\begin{dfn}[Solutions to the Reynolds-defect-equation]
A solution to the Reynolds-defect-equation is a tuple $(u,p,R)$ of smooth functions
\begin{align*}
u\in C^\infty([0,1]\times \R^2;\R^2),p\in C^\infty([0,1]\times \R^2), R\in C^\infty([0,1]\times \R^2;\operatorname{Sym}_{2\times 2})
\end{align*} with
\begin{align*}
\supp (u,p,R )\subset [0,1]\times B_1(0)
\end{align*}
such that
\begin{align*}
\partial_t u + \dv(u\otimes u) + \nabla p &=-\dv \overset{\circ}{R},\\
\dv u &= 0
\end{align*}
is satisfied.
\end{dfn}

\begin{prop}[Main Proposition]\label{prop: main proposition}
There exists a constant $M_0>0$ such that the following holds. Let $e\in C^\infty\left([0,1];\left[\frac{1}{2},1\right]\right)$ be an arbitrary given energy profile.  Choose $\delta>0$  and assume that there exists a smooth solution $(u_0, R_0,p_0)$ to the Reynolds-defect-equation, satisfying
\begin{align}
\frac{3}{4}\delta e(t)&\leq e(t)-\int_{\R^2}|u_0|(t,x)^2\,\dx\leq\frac{5}{4}\delta e(t),\label{eq: assumption energy}\\
40 \|R_0\|_{C_tL^1_x}&< \frac{1}{16}\delta. \label{eq: assumption r}
\end{align}
Then there exists another solution $(u_1, R_1,p_1)$ such that for all $t\in [0,1]$
\begin{enumerate}[(i)]
\item \begin{align*}
\frac{3}{8}\delta e(t)\leq e(t)-\int_{\R^2}|u_1|^2(t,x)\,\dx\leq\frac{5}{8}\delta e(t),
\end{align*}
\item $40\|R_1(t)\|_{L^1(\R^2)}< \frac{1}{32}\delta,$
\item $\|u_1(t)-u_0(t)\|_{L^2(\R^2)}\leq M_0 \delta^\frac{1}{2}$,
\item $\curl(u_1-u_0)\in H^p(\R^2)$ with $\|\curl (u_1-u_0)(t)\|^p_{H^p(\R^2)}\leq \delta$.
\end{enumerate}
\end{prop}
Once Proposition \ref{prop: main proposition} is settled, Theorem \ref{thm:main} is a direct consequence. The method of proof for Theorem \ref{thm:main} is standard in current works on convex integration techniques and we therefore omit it. For a detailed proof of a similar theorem, see \cite{buck2023non}. In order to prepare the proof of Proposition \ref{prop: main proposition}, we will define our building blocks in the next Section \ref{sec: building blocks}. After that, Sections \ref{sec: perturbations} to \ref{sec: proof of main prop} are devoted to the proof of Proposition \ref{prop: main proposition}. In particular, we will define the new solution $u_1$ in Section \ref{sec: perturbations} and $R_1$, $p_1$ in Section \ref{sec: def of R} and fix $M_0$ in Section \ref{sec: proof of main prop}.

\section{The Building Blocks}\label{sec: building blocks}
In this section we define three building blocks $W_k$, $Y_k$, $A_k$. For $p\in (0,1)$, fix $N^\ast=N^\ast(p)\in\N$ maximal such that $N^\ast\leq 2(\frac{1}{p}-1)$. We will need to control the first $N^\ast$ momenta of the vorticity of our perturbations to achieve that, $\curl(u_1-u_0)\in H^p(\R^2)$. For this $N^\ast$, we fix another $N\in\N$ with $2N+3\geq N^\ast.$
Next, fix the vectors $$\xi_1= e_1, \xi_2 = e_2, \xi_3= e_1+e_2, \xi_4 = e_1-e_2$$ in $\R^2$. In the following, we will introduce several parameters that will be fixed in the course of this paper. They will be fixed in the order given by Table \ref{tab: parameters}.\vspace{0,3cm}\\
\begin{table}[H]
\caption{Occuring parameters and their meaning}
\begin{tabular}[H]{l|l}
Parameter & meaning  \\
\hline
$\delta$ & parameters in the Proposition \ref{prop: main proposition} that will ensure convergence\\
$\kappa$ & size of the support of $R_0$ (see \eqref{wishlist: kappa1})\\
$\varepsilon$ &  smoothing of $\rho$ (see \eqref{slow coefficients})\\
$\mu_1$ & concentration (see \eqref{def of base functions}) \\
$\mu_2$ & very high concentration (see \eqref{def of base functions})\\
$\omega$ & phase speed (see \eqref{def of base functions})\\
$\lambda$ & oscillation (see \eqref{def of base functions})
\end{tabular}\label{tab: parameters}
\end{table}
Let $\Phi:\R\rightarrow\R$ be a smooth, odd function with support in $(-\frac{1}{2},\frac{1}{2})$, and $\int\Phi\,\dx = 0$ such that $\varphi := D^{(2N+3)}\Phi$ satisfies $\int\varphi^2\,\dx = 1$. Furthermore, we denote by $\varphi_\mu^k$ the translated functions
\begin{align*}
\varphi_\mu^k(x) = \varphi_\mu\left(x-\frac{k}{16}|\xi_k|^2\right),\quad k=1,2,3,4,
\end{align*}
see also Subsection \ref{subsection: notation}, point (m) for the notation $\varphi_\mu$ and recall the scaling
\begin{align}\label{scaling lambda and mu}
\|D^l\varphi_\mu(\lambda\cdot)\|_{L^s(\mathbb{T})} = \lambda^l\mu^{l+\frac{1}{2}-\frac{1}{s}}, s\in [1,\infty].
\end{align}
The translation will ensure the disjointness of the supports of different building blocks, see  \mbox{Lemma \ref{lemma: supports}} for a precise description of the supports.
Let $\lambda,\omega\gg 1$ with $\lambda\in\N$ and $\mu_2\gg\mu_1\gg 1$ with $\frac{\mu_2}{\mu_1}=\lambda$ to be fixed in \mbox{Section \ref{sec: proof of main prop}.} For $k=1,2,3,4$, let us introduce
\begin{align}
\alpha_k(x) &= \frac{1}{(\lambda\mu_2)^{2N+3}}\varphi^k_{\mu_1}(\lambda x_1)\Phi_{\mu_2}(\lambda x_2),\nonumber\\
v_k(x) &=\partial_2^{2N+3} \alpha_k(x)= \varphi^k_{\mu_1}(\lambda x_1)\varphi_{\mu_2}(\lambda x_2)\nonumber\\
q_k(x) &= \frac{1}{\omega}(\varphi^k_{\mu_1})^2(\lambda x_1)\varphi_{\mu_2}^2(\lambda x_2).\label{def of base functions}
\end{align}

\begin{lemma}\label{lemma: small w est}
It holds
\begin{align*}
 \int_{\tor} v_k^2 \,\dx &= \omega \int_{\tor} q_k\,\dx =1,\\
 \int_{\tor} v_k \,\dx  &= 0.
\end{align*}
For any $s\in [1,\infty]$, we have the estimates
\begin{align*}
\|\partial_1^{l_1}\partial_2^{l_2} \alpha_k\|_{L^s(\tor)} &\leq C(s)\lambda^{l_1+l_2-(2N+3)}\mu_1^{l_1+\frac{1}{2}-\frac{1}{s}}\mu_2^{l_2-(2N+3)+\frac{1}{2}-\frac{1}{s}},\\
\|\partial_1^{l_1}\partial_2^{l_2} v_k\|_{L^s(\tor)} &\leq C(s) \lambda^{l_1+l_2}\mu_1^{l_1+\frac{1}{2}-\frac{1}{s}}\mu_2^{l_2+\frac{1}{2}-\frac{1}{s}},\\
 \|\partial_1^{l_1}\partial_2^{l_2} q_k\|_{L^s(\tor)} &\leq C(s) \omega^{-1}\lambda^{l_1+l_2}\mu_1^{l_1+1-\frac{1}{s}}\mu_2^{l_2+1-\frac{1}{s}}.
\end{align*}
\end{lemma}
\begin{proof}
We have 
\begin{align*}
\int_{\tor} v_k(x) \,\dx =\int_0^1 (\varphi_{\mu_1}^k)^2(\lambda x_1) \,\di\cdot \int_0^1 \varphi^2_{\mu_2}(\lambda x_2) \,\dii =1
\end{align*}
by  \eqref{scaling lambda and mu} and since $\int \varphi^2\, \dx =1$. Since $v_k=\partial_2^{2N+3}\alpha_k$ is a derivative, it holds $\int_{\tor} v_k\,\dx=0$. The estimates can also be proven using \eqref{scaling lambda and mu}.
\end{proof}
For $k=1,2,3,4$ we define the linear rotations
\begin{align}\label{def: Lambda}
\Lambda_k:\R^2&\rightarrow\R^2,
x\mapsto (\xi_k\cdot x, \xi_k^\perp\cdot x).
\end{align}

\begin{lemma}\label{lemma: trafo with rotation}
Let $f\in L^1(\tor)$. Then 
\begin{align*}
\int_{\tor} f(\Lambda_k x) \,\dx = \int_{\tor} f(x)\,\dx
\end{align*}
for all $k=1,2,3, 4$.
\end{lemma}
\begin{proof}
This is immediate for $k=1$ and the same is true for $k=2$ by switching the roles of $x_1$ and $x_2$. For $k=3$, we calculate with the transformation rule by rotating the cube $[-\frac{1}{2},\frac{1}{2}]^2$ by $\Lambda_k$
\begin{align*}
\int_{\tor} f(\Lambda_k x)\,\dx &=  \int_{[-\frac{1}{2},\frac{1}{2}]^2} f(\Lambda_k x)\,\dx\\
&= \frac{1}{|\det D\Lambda_k|}\int_Q f(x) \,\dx 
\end{align*}
where $Q = \Lambda_k([-\frac{1}{2},\frac{1}{2}]^2)$ is the by 90 degress rotated and scaled cube with vertices $\left\lbrace \pm e_1, \pm e_2\right\rbrace$. It is not difficult to see that, by a geometric argument,  it holds $\int_Q f \,\dx= 2\int_{\tor}f \,\dx$ because $f$ is periodic. Since $|\det D\Lambda_k|=2$ for $k=3$, this yields the claim and the same reasoning holds for $k=4$. 
\end{proof}
Our three main building blocks $W_k, Y_k \in C^\infty(\tor;\R^2),$ $A_k\in C^\infty(\tor;\R^{2\times 2})$ depending on $\xi_k,\mu_1,\mu_2,\lambda,\omega$ are now defined as
\begin{align*}
W_k(t,x)&:= v_k\left(\Lambda_k\left(x-\omega t \frac{\xi_k}{|\xi_k|^2}\right)\right) \frac{\xi_k}{|\xi_k|}= v_k\left(\Lambda_k x - \omega t e_1\right)\frac{\xi_k}{|\xi_k|}\\
Y_k(t,x)&:=q_k\left(\Lambda_k\left(x-\omega t \frac{\xi_k}{|\xi_k|^2}\right)\right)\xi_k = q_k\left(\Lambda_k x - \omega t e_1\right)\xi_k,\\
 A_k(t,x) &:= \frac{1}{|\xi_k|^{2N+3}} \alpha_k\left(\Lambda_k\left(x-\omega t \frac{\xi_k}{|\xi_k|^2}\right)\right) \frac{\xi_k}{|\xi_k|}\otimes \frac{\xi_k^\perp}{|\xi_k|}\\
&=\frac{1}{|\xi_k|^{2N+3}} \alpha_k\left(\Lambda_k x-\omega t e_1\right) \frac{\xi_k}{|\xi_k|}\otimes \frac{\xi_k^\perp}{|\xi_k|}
\end{align*}
which means that we first rotate $v_k, q_k$ and $\alpha_k$ and move in time in the direction of $\xi_k$. 

We note that our building blocks are again periodic functions on $\R^2$ with period $1$ in both variables since $\xi_k\in\N^2$.
\begin{prop}[Building blocks]\label{prop: building blocks}
The building blocks are $\lambda$-periodic and  satisfy 
\begin{enumerate}[(i)]
\item $\dv( W_k\otimes W_k) = \partial_t Y_k$,
\item $\int_{\tor} W_k\otimes  W_k (t,x)\,\dx = \frac{\xi_k}{|\xi_k|}\otimes\frac{\xi_k}{|\xi_k|}$,
\item $\int_{\tor} W_k(t,x)\,\dx =0$,
\item  $\int_{\tor} q_k\left(\Lambda_k\left(x-\omega t \frac{\xi_k}{|\xi_k|^2}\right)\right)\,\dx =\frac{1}{\omega}$, $\int_{\tor} Y_k (t,x)\,\dx =\frac{1}{\omega}\xi_k$,
\item $\|W_k(t,\cdot)\|_{L^s(\tor)} = \|v_k\|_{L^s(\tor)}$ for all $s\in [1,\infty]$,
\item  $\|Y_k(t,\cdot)\|_{L^s(\tor)} = \|q_k\|_{L^s(\tor)}$ for all $s\in [1,\infty]$.

\end{enumerate}
\end{prop}
\begin{proof}
For $(i)$, we have by Lemma \ref{lemma: div of matrix} with $f(x) = (\varphi^k_{\mu_1})^2(\lambda (x - \omega t))$ and $g(x)=\varphi_{\mu_2}^2(\lambda x)$
\begin{align*}
\dv( W_k\otimes W_k) & = \dv \left((\varphi^k_{\mu_1})^2(\lambda (\xi_k\cdot x - \omega t))\varphi_{\mu_2}^2(\lambda \xi_k^\perp\cdot x)\frac{\xi_k}{|\xi_k|}\otimes\frac{\xi_k}{|\xi_k|} \right)\\
&= \lambda \left((\varphi^k_{\mu_1})^2\right)'(\lambda (\xi_k\cdot x - \omega t))\varphi_{\mu_2}^2(\lambda\xi_k^\perp\cdot x)\xi_k\\
&=\partial_t Y_k.
\end{align*}
The remaining points follow directly from Lemma \ref{lemma: small w est} and Lemma  \ref{lemma: trafo with rotation}. 
\end{proof}

\begin{prop}\label{prop: building block estimates}
The building blocks satisfy the following estimates uniformly in $t$:
\begin{align*}
\|\nabla^l W_k(t)\|_{L^s(\tor)} &\leq C(s)\lambda^l\mu_1^{\frac{1}{2}-\frac{1}{s}}\mu_2^{l+\frac{1}{2}-\frac{1}{s}},\\
\|\nabla^l Y_k(t)\|_{L^s(\tor2)}  &\leq  C(s) \omega^{-1}\lambda^l\mu_1^{1-\frac{1}{s}}\mu_2^{l+1-\frac{1}{s}},\\
\|\nabla^l A_k(t)\|_{L^s(\tor)}, \|\nabla^l A^T_k(t)\|_{L^s(\tor)}  &\leq C\lambda^{l-(2N+3)}\mu_1^{\frac{1}{2}-\frac{1}{s}}\mu_2^{l-(2N+3)+\frac{1}{2}-\frac{1}{s}},\\
\|\nabla^l \partial_t A_k(t)\|_{L^s(\tor)}, \|\nabla^l \partial_t A^T_k(t)\|_{L^s(\tor)} &\leq C \omega\lambda^{1+l-(2N+3)}\mu_1^{\frac{3}{2}-\frac{1}{s}}\mu_2^{l-(2N+3)+\frac{1}{2}-\frac{1}{s}},\\
\|\nabla^l\dv A_k(t)\|_{L^s(\tor)} &\leq C\lambda^{l-(2N+2)}\mu_1^{\frac{1}{2}-\frac{1}{s}}\mu_2^{l-(2N+2)+\frac{1}{2}-\frac{1}{s}},\\
\|\nabla^l\dv A_k^T(t)\|_{L^s(\tor)} &\leq C\lambda^{l-(2N+2)}\mu_1^{\frac{3}{2}-\frac{1}{s}}\mu_2^{l-(2N+3)+\frac{1}{2}-\frac{1}{s}},\\
\|\nabla^l\dv\dv A_k(t)\|_{L^s(\tor)} &\leq C\lambda^{l-(2N+1)}\mu_1^{\frac{3}{2}-\frac{1}{s}}\mu_2^{l-(2N+2)+\frac{1}{2}-\frac{1}{s}}.
\end{align*}
\end{prop}

\begin{proof}
The estimates follow from Lemma \ref{lemma: trafo with rotation} together with Lemma \ref{lemma: small w est} and we also use that $\mu_2 \gg \mu_1$. For $\dv A_k$, $\dv A_k^T$ and $\dv\dv A_k$ we note that by 
 Lemma \ref{lemma: div of matrix}
\begin{align*}
\dv A_k(t,x) &= \frac{1}{|\xi_k|^{2N+3}}(\partial_2\alpha_k)(\Lambda_k x-\omega t e_1)\xi_k,\\
 \dv\dv A_k(t,x) &= \frac{1}{|\xi_k|^{2N+1}}(\partial_1\partial_2\alpha_k)(\Lambda_k x- \omega t e_1),\\
\dv A_k^T(t,x) &= \frac{1}{|\xi_k|^{2N+3}}(\partial_1\alpha_k)(\Lambda_k x-\omega t e_1)\xi_k^\perp,
\end{align*}
therefore the estimates follow again from Lemma \ref{lemma: trafo with rotation} together with Lemma \ref{lemma: small w est}.
\end{proof}

We will use that for different $k,k'\in \lbrace 1,2,3,4\rbrace$, the building blocks $W_k$, $W_{k'}$ have disjoint support. For the proof of the following Lemma we refer to \cite{buck2023non}, Lemma 4.3.
\begin{lemma}[Disjointness of supports]\label{lemma: supports}
We have 
\begin{align*}
\supp \alpha_k\left(\Lambda_k\left(x-\omega t\frac{\xi_k}{|\xi_k|^2}\right)\right)\subset B_{\frac{1}{\lambda\mu_1}}(0)  +\frac{1}{\lambda}\Lambda_k^{-1}\left(\left(\frac{1}{2},\frac{1}{2}\right)+\frac{k}{16}|\xi_k|^2e_1+\Z^2\right) + \omega t\frac{\xi_k}{|\xi_k|^2},
\end{align*}
and
\begin{align*}
\supp \alpha_k\left(\Lambda_k\left(x-\omega t\frac{\xi_k}{|\xi_k|^2}\right)\right) &=\supp W_k =   \supp Y_k = \supp A_k
\end{align*}
i.e. at time $t$ they are supported in balls of radius $\frac{1}{\lambda\mu_1}$ around the set of points
\begin{align*}
M_k(t)= \left\lbrace\frac{1}{\lambda}\Lambda_k^{-1}\left(\left(\frac{1}{2},\frac{1}{2}\right)+\frac{k}{16}|\xi_k|^2e_1+\Z^2\right) + \omega t\frac{\xi_k}{|\xi_k|^2}\right\rbrace.
\end{align*}
For large enough $\mu_1$ (independent of $\lambda, \mu_2$) it holds
\begin{align*}
\supp \alpha_{k_1}\left(\Lambda_{k_1}\left(x-\omega t\frac{\xi_{k_1}}{|\xi_{k_1}|^2}\right)\right)\cap \supp \alpha_{k_2}\left(\Lambda_{k_2}\left(x-\omega t\frac{\xi_{k_2}}{|\xi_{k_2}|^2}\right)\right)= \emptyset
\end{align*}
for $k_1\neq k_2$.
\end{lemma}

\section{The perturbations}\label{sec: perturbations}
Having defined the building blocks, we start now the proof of Proposition \ref{prop: main proposition}.\\ Let \mbox{$e\in C^\infty([0,1];[\frac{1}{2},1])$,} $\delta>0$ and a smooth solution $(u_0, R_0, p_0)$ to the Reynolds-defect-equation be given with $\supp (u_0, R_0, p_0)\subset [0,1]\times B_1(0)$ and satisfying \eqref{eq: assumption energy} and \eqref{eq: assumption r} as in the statement of Proposition \ref{prop: main proposition}. We have to define $M_0>0$ and a solution $(u_1, R_1, p_1)$ to the Reynolds-defect-equation satisfying $\supp (u_1, R_1, p_1)\subset [0,1]\times B_1(0)$ and properties $(i)$--$(iv)$ in Proposition \ref{prop: main proposition}. In this section, we will define $u_1$, whereas $R_1$, $p_1$ will be defined in Section \ref{sec: def of R}. We will fix the constant $M_0$ in Section \ref{sec: proof of main prop}.\\
Before we can define the perturbations, let us decompose the error $\Rtr$ in the following way. There are smooth functions $\Gamma_k$ with $|\Gamma_k|\leq 1$ such that for any matrix $A$ with $|A-I|<\frac{1}{8}$
\begin{align}\label{eq: decomposition}
A=\sum_k\Gamma_k^2(A)\frac{\xi_k}{|\xi_k|}\otimes\frac{\xi_k}{|\xi_k|},
\end{align}
see Section 5 in \cite{brue2021nonuniqueness}. Recall that $\supp(u_0, R_0, p_0)\subset B_1(0)$.
We choose $\kappa\in (\frac{1}{2},1)$ such that
\begin{align}\label{wishlist: kappa1}
\supp R_0(t,\cdot)\subseteq B_\kappa(0)
\end{align}
for all $t\in[0,1]$. Let $\chi_\kappa\in C^\infty_c(B_1(0))$ be a smooth cutoff with $\chi_\kappa\equiv 1$ on $B_\kappa(0)$. For $\varepsilon>0$ we further define
\begin{align}\label{slow coefficients}
\gamma(t)&= \frac{e(t)(1-\frac{\delta}{2})-\int_{\R^2}|u_0|^2(t,x)\,\dx}{2\|\chi_\kappa\|^2_{L^2(\R^2)}},\nonumber\\
\rho(t,x)&=10\sqrt{\varepsilon^2+|\Rtr(t,x)|^2} + \gamma(t),\nonumber\\
a_k(t,x)&=\chi_\kappa(x)\rho^\frac{1}{2}(t,x)\Gamma_k\left(I+\frac{\Rtr(t,x)}{\rho(t,x)}\right),
\end{align}
noting that the decomposition \eqref{eq: decomposition} exists for $I+\frac{\Rtr}{\rho}$.  For later use, we note that
\begin{align}\label{eq: R coeff}
\chi_\kappa^2(x)\rho(t,x) I +  \Rtr(t,x) = \sum_k a_k^2(t,x)\frac{\xi_k}{|\xi_k|}\otimes\frac{\xi_k}{|\xi_k|}.
\end{align}
We set $$u_{1} = u_0 +w,$$
where we define the perturbation $w$ as $$w = w^p + w^t$$ with:
\begin{align}
w^p &=\sumk \nabla^\perp\curl\Delta^N\left[\dv \left(a_k (A_k+A_k^T)\right)\right],\label{def: wp}\\
w^t &= \sumk \nabla^\perp\curl\Delta^N\left[-\bl^{N+1}_M\left(a_k^2,q_k\left(\Lambda_k \cdot-\omega t e_1\right)-\frac{1}{\omega}\right)\xi_k \right].\label{def: wt}
\end{align}
where $\bl^{N+1}_M$ is defined in Proposition \ref{prop: improved anti laplacian}.
Note that this is well-defined since we have $$\int_{\tor} q_k\left(\Lambda_k\left(x-\omega t \frac{\xi_k}{|\xi_k|^2}\right)\right)\,\dx =\frac{1}{\omega}$$ by Proposition \ref{prop: building blocks}.
As explained in the Introduction, the purpose of writing $w$ as a derivative of order $2N+2$ is that $\curl w$ is a derivative of order $2N+3\geq N^\ast$, which guarantees that all momenta of $\curl w$ up to order $N^\ast$ vanish. This implies $\curl w \in H^p(\R^2)$. We will sometimes use the above definitions \eqref{def: wp} and \eqref{def: wt}, i.e. writing $w$ as $w=\sum_k \nabla^\perp \curl \Delta^N(\dots)$ and use estimates on $A^k$ and its derivatives as well as estimates for $\bl^{N+1}_M$, whereas on other occasions we have to decompose $w= u^p + u^c + u^t+ u^{cc}$ (see Proposition \ref{prop: structure of perturbations} below) and use certain properties of the individual parts. Also note that $\dv w= 0$, being an orthogonal gradient. We recover our building blocks from the previous section from the above definition:
\begin{prop}\label{prop: structure of perturbations}
It holds
\begin{align*}
w = w^p + w^t
\end{align*}
with
\begin{align*}
w^p= u^p+ u^{c}, \quad w^t= u^t  + u^{cc}
\end{align*}
with principal perturbation
\begin{align*}
u^p &= \sum_k a_k W_k,
\end{align*}
time corrector
\begin{align*}
u^t = -\sum_k a_k^2 Y_k,
\end{align*}
and further correctors
\begin{align*}
u^c &=\sum_k\left[\Delta^{N+1}, a_kf_k\right]g_k\xi_k - \sum_k\nabla\dv\left(\Delta^N(a_k\dv A_k)\right)  +\sum_k \nabla^\perp\curl\Delta^N( A_k\cdot\nabla a_k + \dv(a_k A_k^T)),\\
u^{cc} &=   \sum_k \frac{1}{\omega} a_k^2\xi_k- \vl^{N+1}_M\left(a_k^2,q_k\left(\Lambda_k x-\omega t e_1\right)-\frac{1}{\omega}\right)\\
&\hspace{0,3cm} + \sum_k\nabla\dv\Delta^N\left(  \bl^{N+1}_M\left(a_k^2,q_k\left(\Lambda_k x-\omega t e_1\right)-\frac{1}{\omega}\right)\right)\xi_k\\
&=: u^{cc}_I+ u^{cc}_{II} + u^{cc}_{III}
\end{align*}
where we set for simplicity
\begin{align*}
f_k(t,x)&:= \varphi^k_{\mu_1}(\lambda(\xi_k\cdot x-\omega t)),\\
g_k(t,x)&:= \frac{1}{|\xi_k|^{2N+3}(\lambda\mu_2)^{2N+2}}(\Phi')_{\mu_2}(\lambda\xi_k^\perp\cdot x).
\end{align*}
\end{prop}

\begin{proof}
For $w^p$, we see using the elementary identity $\nabla^\perp\curl = \Delta-\nabla\dv$ that
\begin{align}
w^p &=\sum_k\nabla^\perp\curl\Delta^N\dv \left(a_k (A_k+A_k^T)\right)\nonumber\\
&=\sum_k\nabla^\perp\curl\Delta^N( a_k \dv A_k) +\sum_k\nabla^\perp\curl\Delta^N( A_k\cdot\nabla a_k + \dv(a_k A_k^T)) \nonumber\\
&= \sum_k\Delta^{N+1}( a_k \dv A_k)-\sum_k \nabla\dv\left(\Delta^N(a_k\dv A_k)\right)+\sum_k\nabla^\perp\curl\Delta^N( A_k\cdot\nabla a_k + \dv(a_k A_k^T))\label{corrector terms i}
\end{align}
Our main perturbation $u^p$ is "hidden" in the first term, namely when all derivatives from $\Delta^{N+1}$ fall on the part that is concentrated with $\mu_2$. More precisely, we write with $f_k, g_k$ as in the statement of the Proposition and the commutator from \eqref{commutator}
\begin{align}
\sum_k\Delta^{N+1}( a_k \dv A_k) &=\sum_k\Delta^{N+1}\left( a_k \frac{1}{|\xi_k|^{2N+3}}(\partial_2\alpha_k)(\Lambda_k x-\omega t e_1)\xi_k\right)\\
&= \sum_k\Delta^{N+1}\left(a_k f_kg_k\xi_k\right)\nonumber\\
&=\sum_k a_kf_k\Delta^{N+1}(g_k) \xi_k+ \left[\Delta^{N+1},a_kf_k\right] g_k\xi_k\nonumber\\
&= \sum_k a_k \varphi^k_{\mu_1}(\lambda(\xi_k\cdot x-\omega t))\varphi_{\mu_2}(\lambda\xi_k^\perp\cdot x)\frac{\xi_k}{|\xi_k|} + \left[\Delta^{N+1},a_kf_k\right] g_k\xi_k\nonumber\\
&= \sum_k a_k W_k+  \left[\Delta^{N+1},a_kf_k\right] g_k\xi_k = u^p + \sum_k\left[\Delta^{N+1},a_kf_k\right] g_k\xi_k \label{corrector terms ii}.
\end{align}
Combining \eqref{corrector terms i} and \eqref{corrector terms ii}, we have
\begin{align*}
w^p= u^p + u^c
\end{align*}
with
\begin{align*}
u^c &=\sum_k\left[\Delta^{N+1}, a_kf_k\right] - \sum_k\nabla\dv\left(\Delta^N(a_k\dv A_k)\right)  +\sum_k \nabla^\perp\curl\Delta^N( A_k\cdot\nabla a_k + \dv(a_k A_k^T)).
\end{align*}
For $w^t$, we use again the elementary identity $\nabla^\perp\curl = \Delta -\nabla\dv$ and properties of the operator $\bl^{N+1}_M$ that was constructed in Proposition \ref{prop: improved anti laplacian} and obtain with $\vl^{N+1}_M$ from that Lemma:
\begin{align*}
w^t &= -\sum_k\nabla^\perp\curl\Delta^N\bl^{N+1}_M\left(a_k^2,q_k\left(\Lambda_k \cdot-\omega t e_1\right)-\frac{1}{\omega}\right)\xi_k\\
 &\hspace{0,3cm}= -\sum_k \Delta^{N+1}\bl^{N+1}_M\left(a_k^2,q_k\left(\Lambda_k \cdot-\omega t e_1\right)-\frac{1}{\omega}\right)\xi_k\\
&\hspace{0,5cm}  + \sum_k\nabla\dv\Delta^N\left(  \bl^{N+1}_M\left(a_k^2,q_k\left(\Lambda_k \cdot-\omega t e_1\right)-\frac{1}{\omega}\right)\right)\xi_k\\
&\hspace{0,3cm} =- \sum_k a_k^2\left(q_k\left(\Lambda_k x-\omega t e_1\right)-\frac{1}{\omega}\right)\xi_k- \vl^{N+1}_M\left(a_k^2,q_k\left(\Lambda_k \cdot-\omega t e_1\right)-\frac{1}{\omega}\right)\xi_k\\
&\hspace{0,5cm} + \sum_k\nabla\dv\Delta^N\left(  \bl^{N+1}_M\left(a_k^2,q_k\left(\Lambda_k \cdot-\omega t e_1\right)-\frac{1}{\omega}\right)\right)\xi_k\\
&\hspace{0,3cm}= -\sum_k a_k^2 Y_k + \sum_k \frac{1}{\omega} a_k^2\xi_k- \vl^{N+1}_M\left(a_k^2,q_k\left(\Lambda_k \cdot-\omega t e_1\right)-\frac{1}{\omega}\right)\xi_k\\
&\hspace{0,5cm} + \sum_k\nabla\dv\Delta^N\left(  \bl^{N+1}_M\left(a_k^2,q_k\left(\Lambda_k \cdot-\omega t e_1\right)-\frac{1}{\omega}\right)\right)\xi_k\\
&\hspace{0,3cm}= u^t + u^{cc}.
\end{align*}
This yields
\begin{align*}
w^t=  u^t  + u^{cc}.
\end{align*}
\end{proof}

\section{Estimates of the perturbations}\label{sec: estimates perturbations}\label{sec. estimates of perturbations}
In this section, we provide the necessary estimates on the perturbations. We start with a preliminary estimate for the coefficients $a_k$ and then estimate the individual parts of the perturbations separately. After that, we obtain an estimate on the energy increment and conclude the section by fixing the parameter $\varepsilon$.
\begin{lemma}[Preliminary estimates]
It holds
\begin{align}
\|a_k\|_{C^l(\R^2\times [0,1])}&\leq C(R_0, u_0, e,\delta,\kappa,\varepsilon,l),\label{est ak}
\end{align}
and
\begin{align}\label{eq: a_k in L2}
\|a_k(t,\cdot)\|_{L^2(\R^2)}\leq \sqrt{10\pi}\left(\varepsilon^\frac{1}{2} + \delta^\frac{1}{2}\right)
\end{align}
uniformly in $t$.
\end{lemma}
\begin{proof}
Since $R_0, u_0, e$ and the functions $\Gamma_k$ are smooth, it is clear that also $a_k$ is smooth and that its derivatives can be estimated by the entities $R_0, u_0, e,\delta,\kappa,\varepsilon$. For the second part, we calculate using $|\Gamma_k|\leq 1$
\begin{align*}
\int_{\R^2} a_k^2(t,x) \,\dx &= \int_{\R^2}\chi_\kappa^2(x)\rho(t,x) \Gamma_k^2\left(I+\frac{\Rtr(t,x)}{\rho(t,x)}\right)\,\dx \leq \int_{B_1(0)} 10\sqrt{\varepsilon^2+|\Rtr(t)|^2(t,x)}+\gamma(t)  \,\dx\\
& \leq 10\pi\varepsilon + 20\|R_0\|_{C_tL^1_x} + 10 \pi  \gamma(t)\\
&\leq 10\pi\varepsilon + 20\|R_0\|_{C_tL^1_x} +  \frac{5\pi}{\|\chi_\kappa\|^2_{L^2(\R^2)}}\left(e(t)(1-\frac{\delta}{2}) - \int_{\R^2} |u_0|^2(t,x)\,\dx\right).
\end{align*}
We have
\begin{align*}
0\leq e(t)\left(1-\frac{\delta}{2}\right) - \int_{\R^2} |u_0|^2(t,x)\,\dx  \leq \frac{3}{4}\delta,
\end{align*}
by assumption. Using \eqref{eq: assumption energy}, $\kappa>\frac{1}{2}$ and 
$$\frac{5\pi}{\|\chi_\kappa\|^2_{L^2(\R^2)}}\leq  \frac{5}{\kappa^2}\leq 20,$$ 
and also that $20\|R_0\|_{C_tL^1_x}\leq \delta$ by assumption,
we obtain
\begin{align}\label{eq: a_k in L2 quadr}
\int_{\R^2} a_k^2(t,x) \,\dx \leq 10\pi\varepsilon   + 16\delta.
\end{align}
From this \eqref{eq: a_k in L2} follows.
\end{proof}

\begin{lemma}[Estimate of the principal perturbation]\label{lemma: estimate up}
It holds
\begin{align*}
\|u^p(t)\|_{L^s(\R^2)}&\leq C(R_0,u_0,e,\delta,\kappa, \varepsilon)\mu_1^{\frac{1}{2}-\frac{1}{s}}\mu_2^{\frac{1}{2}-\frac{1}{s}}\text{ for } s\in [1,\infty]
\end{align*}
and for $s=2$  more refined
\begin{align}
\|u^p(t)\|_{L^2(\R^2)}&\leq 4\sqrt{10\pi}\left(\varepsilon^\frac{1}{2} + \delta^\frac{1}{2}\right) +\frac{ C(R_0,u_0,e,\delta,\varepsilon)}{\lambda^\frac{1}{2}}\label{wishlist vareps 2}
\end{align}
uniformly in $t$. 
\end{lemma}
\begin{proof}
We have defined $u^p$ as $u^p = \sum_k a_k W_k.$ For the first estimate, we use  the fact that $u^p$ is supported in $[0,1]\times B_1(0)$ (because $a_k$ is supported in $[0,1]\times B_1(0)$), and get
\begin{align*}
\|u^p(t)\|_{L^s(\R^2)} \leq \sum_k \|a_k\|_{C(\R^2\times [0,1])}\|W_k(t)\|_{L^s(B_1(0))} \leq 4\sum_k \|a_k\|_{C(\R^2\times [0,1])}\|W_k(t)\|_{L^s(\tor)},
\end{align*}
so that the first estimate follows from \eqref{est ak} and Proposition \ref{prop: building blocks}.
 For the second estimate, we use Proposition \ref{lemma: improved hoelder}, noting again that $\supp a_k(t,\cdot)\subset B_1(0)$, Lemma \ref{lemma: small w est}, Proposition \ref{prop: building blocks} and \eqref{eq: a_k in L2}
\begin{align*}
\|\sumk a_k(t,\cdot) W_k(t,\cdot)\|_{L^2(\R^2)} &\leq \sum_k\left(\|a_k(t,\cdot)\|_{L^2(\R^2)} \|W_k(t)\|_{L^2(\tor)} + \frac{C}{\lambda^\frac{1}{2}}\|a_k(t,\cdot)\|_{C^1(\R^2)}\|W_k(t)\|_{L^2(\tor)}\right) \\
&\leq 4\sqrt{10\pi}\left(\varepsilon^\frac{1}{2} + \delta^\frac{1}{2}\right) + \frac{C(R_0, u_0, e,\delta, \kappa,\varepsilon)}{\lambda^\frac{1}{2}}.\qedhere
\end{align*}.
\end{proof}

\begin{lemma}[Estimate of $u^c$]\label{lemma: estimate correctors i}
We have
\begin{align*}
\|u^c(t)\|_{L^2(\R^2)}&\leq C(R_0,u_0,e,\delta,\kappa,\varepsilon)\frac{\mu_1}{\mu_2}
\end{align*}
uniformly in $t$.
\end{lemma}
\begin{proof}
We remember that
\begin{align*}
u^c &=\sum_k\left[\Delta^{N+1}, a_kf_k\right] g_k\xi_k - \sum_k\nabla\dv\left(\Delta^N(a_k\dv A_k)\right)\\
&\hspace{0,3cm}  +\sum_k \nabla^\perp\curl\Delta^N( A_k\cdot\nabla a_k) + \sum_k \nabla^\perp\curl\Delta^N\dv(a_k A_k^T)\\
&=: I+ II+ III+ IV.
\end{align*}
For $I$, we remember that  
\begin{align*}
f_k(t,x)&=\varphi^k_{\mu_1}(\lambda(\xi_k\cdot x-\omega t)),\\
g_k(t,x)&= \frac{1}{|\xi_k|^{2N+3}(\lambda\mu_2)^{2N+2}}(\Phi')_{\mu_2}(\lambda\xi_k^\perp\cdot x).
\end{align*}
We use formula \eqref{commutator sum} for the commutator, so that by first applying Lemma \ref{lemma: trafo with rotation}, the estimate \eqref{est ak} and $\supp_x a_k\subseteq B_1(0)$ and the scaling property \eqref{scaling lambda and mu} it holds
\begin{align*}
\|I\|_{L^2(\R^2)}&\leq \sum_k \|\left[\Delta^{N+1},a_k f_k\right]g_k\|_{L^2(\R^2)}\\
& \leq C(R_0, u_0, e,\delta,\kappa,\varepsilon)\sum_k \sum_{\substack{m_1+m_2+m_3=N+1\\m_3\leq N}}\sum_{|\beta|=m_1} \|\partial^\beta\Delta^{m_2}f_k \partial^\beta\Delta^{m_3}g_k \|_{L^2(\tor)}\\
&\leq \frac{C(R_0, u_0, e,\delta,\kappa,\varepsilon)}{(\lambda\mu_2)^{2N+2}}\sum_k \sum_{\substack{m_1+m_2+m_3=N+1\\m_3\leq N}} \|\varphi^k_{\mu_1}(\lambda\cdot)\|_{W^{m_1+2m_2,2}(\mathbb{T})} \|(\Phi')_{\mu_2}(\lambda\cdot) \|_{W^{m_1+2m_3,2}(\mathbb{T})}\\
&\leq \frac{C(R_0, u_0, e,\delta,\kappa,\varepsilon)}{(\lambda\mu_2)^{2N+2}}\sum_k \sum_{\substack{m_1+m_2+m_3=N+1\\m_3\leq N}} \lambda^{2N+2}\mu_1^{m_1+2m_2}\mu_2^{m_1+2m_3}.
\end{align*}
In the third line we first used  Lemma \ref{lemma: trafo with rotation} and then Fubinis theorem to separate the integrals.
  Remembering that $\mu_2\gg \mu_1$, this yields that the summands with the worst estimates are those with the most possible appearances of $\mu_2$. Since $m_1+2m_3\leq 2N+1$, this is maximal $2N+1$, which is achieved only for $(m_1,m_2,m_3)=(1,0,N)$. Hence,
\begin{align*}
\|I\|_{L^2(\R^2)}\leq C(R_0, u_0, e,\delta,\kappa,\varepsilon) \frac{\mu_1}{\mu_2}.
\end{align*}
For $II$, we interchange the order of derivatives  to compute
\begin{align*}
\nabla\dv\left(\Delta^N(a_k\dv A_k)\right) &= \nabla\Delta^N(a_k\dv\dv A)  + \nabla\Delta^N(\dv A_k \cdot \nabla a_k).
\end{align*}
Therefore, it follows by inequality \eqref{est ak}  with $\supp_x a_k\subseteq B_1(0)$ and Proposition \ref{prop: building block estimates} 
\begin{align*}
\|II\|_{L^2(R^2)} &\leq C(R_0, u_0, e,\delta,\kappa,\varepsilon) \left(\|\dv\dv A_k\|_{W^{2N+1,2}(\tor)} + \|A_k\|_{W^{2N+2,2}(\tor)}\right)\\
&\leq C(R_0, u_0, e,\delta,\kappa,\varepsilon)\left(\frac{\mu_1}{\mu_2} + \frac{1}{\lambda\mu_2} \right)\\
&\leq C(R_0, u_0, e,\delta,\kappa,\varepsilon) \frac{\mu_1}{\mu_2}.
\end{align*}
For $(III)$, we use again Proposition \ref{prop: building block estimates}  and \eqref{est ak} and get
\begin{align*}
\|III\|_{L^2(\R^2}&\leq C(R_0, u_0, e,\delta,\kappa,\varepsilon)\|A_k\|_{W^{2N+2,2}(\tor)}\\
&\leq C(R_0, u_0, e,\delta,\kappa,\varepsilon)\frac{1}{\lambda\mu_2}.
\end{align*}
Another use of   Proposition \ref{prop: building block estimates}  and \eqref{est ak} yields for $(IV)$
\begin{align*}
\|IV\|_{L^2(\R^2)}&\leq C(R_0, u_0, e,\delta,\kappa,\varepsilon)\|\dv A_k^T\|_{W^{2N+2,2}(\tor)}  \\
&\leq  C(R_0, u_0, e,\delta,\kappa,\varepsilon)\frac{\mu_1}{\mu_2}.
\end{align*}
\end{proof}

For the time corrector $u^t$ and its corrector terms $u^{cc}$, it is convenient to already estimate the curl in $L^\infty$, which we do in the next two Lemmas. Also we will need a few more estimates, so we get that done now.
\begin{lemma}[Estimate of the time corrector]\label{lemma: estimate time corrector}
It holds 
\begin{align*}
\|u^t(t)\|_{L^1(\R^2)}&\leq C(R_0,u_0,e,\delta,\kappa,\varepsilon)\frac{1}{\omega},\\
\|u^t(t)\|_{L^2(\R^2)}&\leq C(R_0,u_0,e,\delta,\kappa,\varepsilon)\frac{\mu_1^\frac{1}{2}\mu_2^\frac{1}{2}}{\omega},\\
\|\curl u^t(t)\|_{L^\infty(\R^2)}&\leq C(R_0,u_0,e,\delta,\kappa,\varepsilon)\frac{\lambda\mu_1\mu_2^2}{\omega}
\end{align*}
uniformly in $t$. 
\end{lemma}
\begin{proof}
This follows from Proposition \ref{prop: building block estimates} and \eqref{est ak}.
\end{proof}

For the estimate of the corrector $u^{cc}$, recall that  $u^{cc}$ was defined as 
\begin{align*}
u^{cc} &=   \sum_k \frac{1}{\omega} a_k^2\xi_k- \vl^{N+1}_M\left(a_k^2,q_k\left(\Lambda_k \cdot-\omega t e_1\right)-\frac{1}{\omega}\right)\xi_k\\
&\hspace{0,3cm} + \sum_k\nabla\dv\Delta^N\left(  \bl^{N+1}_M\left(a_k^2,q_k\left(\Lambda_k \cdot-\omega t e_1\right)-\frac{1}{\omega}\right)\right)\xi_k\\
&= u^{cc}_I + u^{cc}_{II}+ u^{cc}_{III}.
\end{align*} 
\begin{lemma}[Estimate of $u^{cc}$]\label{lemma: estimate correctors ii}
We have
\begin{align*}
\|u^{cc}(t)\|_{L^2(\R^2)}&\leq  C(R_0, u_0, e,\delta,\kappa,\varepsilon,N,M)\frac{\mu_1^{\frac{1}{2}}\mu_2^{\frac{1}{2}}}{\omega},\\
\|\curl u^{cc}(t)\|_{L^\infty(\R^2)}&\leq C(R_0, u_0, e,\delta,\kappa,\varepsilon,N,M)\left(\omega^{-1}\lambda^{-M+1}\mu_1\mu_2^2 + \frac{1}{\omega}\right)
\end{align*}
uniformly in $t$, and for the part $u^{cc}_{II}$ we have
\begin{align*}
\|\partial_t u^{cc}_{II}(t)\|_{L^1(\R^2)} &\leq  C(R_0, u_0, e,\delta,\kappa,\varepsilon,N,M)\left(\lambda^{-M+1}\mu_1+\frac{1}{\omega}\right)
\end{align*}
uniformly in $t$.
\end{lemma}

\begin{proof}
The expression $u^{cc}_I$ is simply estimated by \eqref{est ak} with
\begin{align*}
\|u^{cc}_I\|_{W^{m,s}(\R^2)}\leq C(R_0, u_0, e,\delta,\kappa,\varepsilon)\frac{1}{\omega}.
\end{align*}
For $u^{cc}_{II}$, we first note that $q_k$ is $\frac{1}{\lambda}$-periodic, i.e. $q_k =\tilde{q}_k(\lambda\cdot)$ for $\tilde{q}_k(x) = \frac{1}{\omega}(\varphi^k_{\mu_1})^2(x_1)\varphi_{\mu_2}^2(x_2)$, see \eqref{def of base functions}. We use \eqref{est: bl2}, \eqref{est ak}, Lemma \ref{lemma: trafo with rotation} and Lemma \ref{lemma: small w est} (that also works for $\tilde{q}_k$ with $\tilde{\lambda}=1$)
\begin{align*}
\|u^{cc}_{II}\|_{W^{m,s}(\R^2)}&= \|\vl^{N+1}_M\left(a_k^2,q_k\left(\Lambda_k \cdot-\omega t e_1\right)-\frac{1}{\omega}\right)\|_{W^{m,s}(\R^2)}\\
  &\hspace{0,3cm}\leq C(N,M) \lambda^{-M+m} \|a_k^2\|_{C^{2(N+1)(M+1)+m}(\R^2)}\|\tilde{q}_k\left(\Lambda_k \cdot-\omega t e_1\right)-\frac{1}{\omega}\|_{W^{m,s}(\tor)}\\
  &\hspace{0,3cm}\leq C(N,M) \lambda^{-M+m} \|a_k^2\|_{C^{2(N+1)(M+1)+m}(\R^2)}\left(\|\tilde{q}_k\|_{W^{m,s}(\tor)} + \frac{1}{\omega}\right)\\
   &\hspace{0,3cm}\leq  C(R_0, u_0, e,\delta,\varepsilon,\kappa, N,M)\left(\omega^{-1}\lambda^{-M+m}\mu_1^{1-\frac{1}{s}}\mu_2^{m+1-\frac{1}{s}} + \frac{1}{\omega}\right).
 \end{align*}
Since we also want to estimate the time derivative of $u^{cc}_{II}$, we estimate  $\partial_t \vl^{N+1}_M$ in the same way using \eqref{properties b v} and $\partial_t \left(q_k(\Lambda_k x- \omega t e_1)\right) = \omega (\partial_1q_k)(\Lambda_k x- \omega t e_1)$
 \begin{align*}
 &\|\partial_t \vl^{N+1}_M\left(a_k^2,q_k\left(\Lambda_k x-\omega t e_1\right)-\frac{1}{\omega}\right)\|_{W^{m,s}(\R^2)}\\
   &\hspace{0,3cm}\leq  C(R_0, u_0, e,\delta,\varepsilon,\kappa,N,M)\left(\lambda^{-M+m+1}\mu_1^{2-\frac{1}{s}}\mu_2^{m+1-\frac{1}{s}}+\frac{1}{\omega}\right).
 \end{align*}
Summarising what we have so far,
\begin{align*}
\|u^{cc}_I\|_{L^2(\R^2)}&\leq  C(R_0, u_0, e,\delta,\kappa,\varepsilon) \frac{1}{\omega}\leq C(R_0, u_0, e,\delta,\kappa,\varepsilon)\frac{\mu_1^\frac{1}{2}\mu_2^\frac{1}{2}}{\omega},\\
\| u^{cc}_{II}\|_{L^2(\R^2)}&\leq C(R_0, u_0, e,\delta,\kappa,\varepsilon,N,M)\left(\omega^{-1}\lambda^{-M}\mu_1^{\frac{1}{2}}\mu_2^{\frac{1}{2}} + \frac{1}{\omega}\right)\\
&\leq C(R_0, u_0, e,\delta,\kappa,\varepsilon,N,M)\frac{\mu_1^\frac{1}{2}\mu_2^\frac{1}{2}}{\omega},\\
\|\curl(u^{cc}_I)\|_{L^\infty(\R^2)}&\leq C(R_0, u_0, e,\delta,\kappa,\varepsilon) \frac{1}{\omega},\\
\|\curl(u^{cc}_{II})\|_{L^\infty(\R^2)}&\leq C(R_0, u_0, e,\delta, \kappa,\varepsilon,N,M)\left(\omega^{-1}\lambda^{-M+1}\mu_1\mu_2^2 + \frac{1}{\omega}\right),\\
\|\partial_t u^{cc}_{II}\|_{L^1(\R^2)} &\leq  C(R_0, u_0, e,\delta,\varepsilon,\kappa, N,M)\left(\lambda^{-M+1}\mu_1+\frac{1}{\omega}\right).
\end{align*} 
 For $u^{cc}_{III}$, let us first note that $\curl u^{cc}_{III} = 0$, since this is true for any gradient. Therefore it suffices to estimate $u^{cc}_{III}$ in $L^2$. We use \eqref{est: bl1}, the classical Calderón-Zygmund estimate, \eqref{est ak}, Lemma \ref{lemma: trafo with rotation} and Lemma \ref{lemma: small w est} for $\tilde{q}_k$ with $\tilde{\lambda}=1$ to obtain
 \begin{align*}
 \|u^{cc}_{III}\|_{L^2(\R^2)}&\leq \|\bl^{N+1}_M(a_k^2,q_k\left(\Lambda_k \cdot-\omega t e_1\right)-\frac{1}{\omega})\|_{W^{2N+2,2}(\R^2)}\\
 &\leq C(N,M) \|a_k^2\|_{C^{2(N+1)M+2N+2}(\R^2)}\|\Delta^{-(N+1)}\left(\tilde{q}_k\left(\Lambda_k \cdot-\omega t e_1\right)-\frac{1}{\omega}\right)\|_{W^{2N+2,2}(\R^2)}\\
 &\leq C(N,M) \|a_k^2\|_{C^{2(N+1)M+2N+2}(\R^2)}\|\tilde{q}_k\left(\Lambda_k \cdot-\omega t e_1\right)-\frac{1}{\omega}\|_{L^2(\R^2)}\\
 &\leq C(R_0, u_0, e,\delta, \kappa,\varepsilon,N,M)\omega^{-1}\mu_1^{\frac{1}{2}}\mu_2^{\frac{1}{2}}.
 \end{align*}
\end{proof}

\begin{lemma}[Estimate of the energy increment]\label{lemma: energy inc}
We have 
\begin{align}
\left|e(t) \left(1-\frac{\delta}{2}\right) - \int_{\R^2} |u_1|^2(t,x)\,\dx\right|&\leq  20 \pi  \varepsilon +  40 \|R_0\|_{C_tL^1_x}\nonumber\\
&\hspace{0,3cm} +  C(R_0, u_0, e,\delta, \kappa,\varepsilon,N,M)\left( \frac{\mu_1^{\frac{1}{2}}\mu_2^{\frac{1}{2}}}{\omega} +\frac{\mu_1}{\mu_2} + \mu_1^{-\frac{1}{2}}\mu_2^{-\frac{1}{2}}+\frac{1}{\lambda} \right).\label{wishlist eps sigma}
\end{align}
\end{lemma}

\begin{proof}
This is proven by standard techniques, see for example \cite{buck2023non}, Lemma 6.6.
\end{proof}
At this point, we fix $\varepsilon$ and choose this parameter so small such that
\begin{align*}
 40 \|R_0\|_{C_tL^1_x}+ 20 \pi  \varepsilon &<\frac{1}{16}\delta,\\
4\sqrt{10\pi} (\varepsilon^\frac{1}{2}+\delta^{\frac{1}{2}})&\leq 40\delta^{\frac{1}{2}},
\end{align*}
which is possible since $ 40 \|R_0\|_{C_tL^1_x}<\frac{1}{16}\delta$ by assumption. Therefore \eqref{wishlist vareps 2} becomes
\begin{align}\label{est up with eps}
\|u^p(t)\|_{L^2(\R^2)}&\leq 40 \delta^\frac{1}{2} +\frac{ C(R_0,u_0,e,\delta,\kappa,\varepsilon)}{\lambda^\frac{1}{2}}
\end{align}
and \eqref{wishlist eps sigma} reduces to
\begin{align}\label{energy increment}
&\left|e(t) \left(1-\frac{\delta}{2}\right) - \int_{\R^2} |u_1|^2(t,x)\,\dx\right|\nonumber\\
&\hspace{0,3cm}< \frac{1}{16}\delta  +  C(R_0, u_0, e,\delta,\kappa,\varepsilon,N,M)\left( \frac{\mu_1^{\frac{1}{2}}\mu_2^{\frac{1}{2}}}{\omega}+\frac{\mu_1}{\mu_2} + \mu_1^{-\frac{1}{2}}\mu_2^{-\frac{1}{2}}+\frac{1}{\lambda} \right)\nonumber\\
&\hspace{0,3cm}\leq \frac{1}{8}\delta e(t) +C(R_0, u_0, e,\delta,\kappa,\varepsilon,N,M)\left( \frac{\mu_1^{\frac{1}{2}}\mu_2^{\frac{1}{2}}}{\omega} +\frac{\mu_1}{\mu_2} + \mu_1^{-\frac{1}{2}}\mu_2^{-\frac{1}{2}}+\frac{1}{\lambda} \right),
\end{align}
where we used that $e(t)\geq\frac{1}{2}$ by assumption.

\section{Estimates of the curl in Hardy space}\label{sec: curl estimates}
In the following Lemmas, we prove that the curls of the perturbations are in the real Hardy space $H^p(\R^2)$ and estimate their Hardy space seminorms in terms of $\lambda, \mu_1$ and $\mu_2$. We will use  \mbox{Remark \ref{rem: hardy atoms};} therefore, we decompose the perturbations into finitely many functions that are supported on disjoint, very small balls of radius $\frac{1}{\lambda\mu_1}$.
\begin{lemma}[Curl of $w^p$]\label{lemma: curl up}
It holds $\curl (w^p)(t)\in H^p(\R^2)$ and
\begin{align*}
\|\curl (w^p)(t)\|_{H^p(\R^2)}\leq C(R_0,u_0,e,\delta,\varepsilon) \lambda\mu_1^{\frac{1}{2}-\frac{2}{p}}\mu_2^\frac{3}{2}\text{ for all } t\in[0,1].
\end{align*}
\end{lemma}

\begin{proof}
By definition of the perturbations and Proposition \ref{prop: structure of perturbations},
\begin{align*}
\curl(w^p)=\sum_k \curl\left(\nabla^\perp\curl\Delta^N\dv \left(a_k (A_k+A_k^T)\right)\right) .
\end{align*}
Since this is a derivative of order $2N+4\geq N^\ast$, we see that all momenta of this expression up to order $N^\ast$ vanish. By the choice of $N^\ast$, this together with the fact that $\supp_x \curl(w^p)(t,\cdot)\subset B_1(0)$ (because $\supp_x a_k\subset B_1(0)$) implies by Remark \ref{rem: hardy atoms} that $\curl(w^p)(t)\in H^p(\R^d)$.
By definition of $A_k$, 
\begin{align*}
\supp_x A_k(t)= \supp_x A_k^T(t) = \supp_x \alpha_k\left(\Lambda_k\left(\cdot-\omega t \frac{\xi_k}{|\xi_k|^2}\right)\right)  .
\end{align*}
Therefore, as seen in  Lemma \ref{lemma: supports}, for a fixed time $t$, the function $(w^p)(t)$ is supported in small, disjoint balls of radius $\frac{1}{\lambda\mu_1}$ around the points in the  finite set $$M(t)= \left\lbrace\frac{1}{\lambda}\Lambda_k^{-1}\left(\left(\frac{1}{2},\frac{1}{2}\right)+\frac{k}{16}|\xi_k|^2e_1+\Z^2\right) + \omega t\frac{\xi_k}{|\xi_k|^2}, k=1,2,3,4\right\rbrace\cap B_1(0).$$ For the rest of the proof, let us fix a time $t$. Let us abbreviate $B_{x_0} = B_{\frac{1}{\lambda\mu_1}}(x_0)$ for $x_0\in M(t)$, and let us decompose $w^p(t)$ as
\begin{align*}
w^p(t,x)= \sum_{x_0\in M(t)} \theta_{x_0}(t,x)
\end{align*}
where
\begin{align*}
\theta_{x_0}(t,x) =\mathbb{1}_{B(x_0)}(x)w^p(t,x).
\end{align*}
The function $\theta_{x_0}(t,\cdot)$ is smooth, it is still a derivative of order $2N+3\geq N^\ast$ and has compact support. Hence $\curl\theta_{x_0}(t,\cdot)\in H^p(\R^2)$. We estimate the $H^p$-seminorm for each $\curl\theta_{x_0}(t,\cdot).$ We have 
\begin{align*}
\curl \theta_{x_0}(t,\cdot) &= \mathbb{1}_{B(x_0)}\curl  (w^p)(t,\cdot) = \mathbb{1}_{B(x_0)}\sum_k \curl\left(\nabla^\perp\curl\Delta^N\dv \left(a_k (A_k+A_k^T)\right)\right)(t,\cdot).
\end{align*}
 As already said, $\theta_{x_0}(t,\cdot)$ is supported on one ball of measure $\frac{C}{\lambda\mu_1}$. We want to use estimate \eqref{est scaled atom}: By Proposition \ref{prop: building block estimates} and \eqref{est ak}, it holds
\begin{align*}
\|\curl\theta_{x_0}(t)\|_{L^\infty(\R^2)}\leq  C(R_0, u_0, e,\delta, \kappa,\varepsilon)\sumk\|A_k+A_k^T(t)\|_{W^{2N+4,\infty}(\R^2)}\leq C(R_0,u_0,e,\delta,\varepsilon)\lambda \mu_{1}^{\frac{1}{2}}\mu_2^{\frac{3}{2}}.
\end{align*}
This gives us by Remark \ref{rem: hardy atoms}
\begin{align*}
\|\curl\theta_{x_0}(t)\|_{H^p(\R^2)}\leq C(R_0,u_0,e,\delta,\kappa,\varepsilon) \lambda^{1-\frac{2}{p}}\mu_1^{\frac{1}{2}-\frac{2}{p}}\mu_2^{\frac{3}{2}}.
\end{align*}
Since $|M(t)|$ is of order $\lambda^2$, $\curl (w^p)$ is made up of $\approx\lambda^2$- many functions $\curl\theta_{x_0}$, and we obtain
\begin{align*}
\|\curl (w^p)(t)\|^p_{H^p(\R^2)}&\leq \sum_{x_0\in M(t)}\|\curl\theta_{x_0}(t)\|^p_{H^p(\R^2)} \leq C(R_0,u_0,e,\delta,\kappa,\varepsilon) \lambda^2 \lambda^{p-2} \mu_1^{\frac{p}{2}-2}\mu_2^\frac{3p}{2}\\ &=  C(R_0,u_0,e,\delta,\kappa,\varepsilon) \lambda^{p} \mu_1^{\frac{p}{2}-2}\mu_2^\frac{3p}{2} .\qedhere
\end{align*}
\end{proof}

\begin{lemma}[Curl of $w^t$]\label{lemma: curl ut}
It holds $\curl (w^t)(t)\in H^p(\R^2)$ and
\begin{align*}
\|\curl w^{t}(t)\|^p_{H^p(\R^2)}\leq C(R_0,u_0,e,\delta,\kappa,\varepsilon, N,M)&\left( \omega^{-p}\lambda^p\mu_1^{p-2}\mu_2^{2p} + \omega^{-p}\lambda^{p(-M+1)}\mu_1^p\mu_2^{2p}\right.\\
&\hspace{0,3cm}\left. + \lambda^2 |\log(\lambda\mu_1)| \omega^{-p} \right)
\end{align*}
for all $t$.
\end{lemma}

\begin{proof}
We remember that by Proposition \ref{prop: structure of perturbations},
\begin{align*}
\curl(w^t)=-\sum_k\curl\left(\nabla^\perp\curl\Delta^N\bl^{N+1}_M\left(a_k^2,q_k\left(\Lambda_k \cdot-\omega t e_1\right)-\frac{1}{\omega}\right)\xi_k\right) =\curl(u^t) +\curl(u^{cc}).
\end{align*}
Since this is a derivative of order $2N+3\geq N^\ast$, we see that all momenta of $\curl(w^t)$ up to order $N^\ast$ vanish (but this might not be true for $\curl(u^t)$ and $\curl(u^{cc})$). By the choice of $N^\ast$, this together with the fact that $\supp_x \curl(w^t)(t,\cdot)\subset B_1(0)$ (because $\supp_x a_k\subset B_1(0)$ and Proposition \ref{prop: improved anti laplacian}) implies by Remark \ref{rem: hardy atoms} that $\curl(w^t)(t)\in H^p(\R^d)$. We want to use again \eqref{est scaled atom}, but here we need to be extra careful with the supports. By Lemma \ref{lemma: supports}, $u^t(t)$ is supported in small, disjoint balls of radius $\frac{1}{\lambda\mu_1}$ around the points in the  finite set $$M(t)= \left\lbrace\frac{1}{\lambda}\Lambda_k^{-1}\left(\left(\frac{1}{2},\frac{1}{2}\right)+\frac{k}{16}|\xi_k|^2e_1+\Z^2\right) + \omega t\frac{\xi_k}{|\xi_k|^2}, k=1,2,3,4\right\rbrace\cap B_1(0).$$ Let us fix a time $t$ and abbreviate again  $B_{x_0} = B_{\frac{1}{\lambda\mu_1}}(x_0)$ for $x_0\in M(t)$, we decompose $u^t(t)$ as
\begin{align*}
u^t(t,x)& = \sum_{x_0\in M(t)} \tilde{\theta}_{x_0}(t,x),\\
\tilde{\theta}_{x_0}(t,x) &=\mathbb{1}_{B(x_0)}(x) u^t(t,x),
\end{align*}
just as we did for $w^p(t)$ in the previous Lemma and note that $\tilde{\theta}_{x_0}(t,\cdot)$ is smooth. On the other hand, $u^{cc}$ is only supported in $B_1(0)$, but not on some smaller balls. We will now apply \mbox{Lemma \ref{Hardy for sum}} with $R=1$, more precisely \eqref{est hardy sum small eps}, to the sum $ \sum_{x_0\in M(t)} \curl\tilde{\theta}_{x_0}(t,\cdot) + \curl u^{cc}(t,\cdot)$, i.e. 
\begin{align}
\| \sum_{x_0\in M(t)} \curl\tilde{\theta}_{x_0}(t,\cdot) + \curl u^{cc}(t,\cdot)\|_{H^p(\R^2)}^p&\leq C(p) \left(\sum_{x_0\in M(t)}\left(\frac{1}{\lambda\mu_1}\right)^2 \|\curl\tilde{\theta}_{x_0}(t,\cdot)\|^p_{L^\infty(\R^2)}\right.\nonumber\\
&\hspace{1cm} +  \|\curl u^{cc}(t,\cdot)\|^p_{L^\infty(\R^2)}\nonumber\\
&\hspace{0,3cm}+ \sum_{x_0\in M(t)} |\log(\lambda\mu_1)| \max_{|\beta|\leq N^\ast}\left|\int x^\beta\curl \tilde{\theta}_{x_0}(t,x) \,\dx\right|^p\nonumber\\
&\hspace{0,3cm}\left.  +\max_{|\beta|\leq N^\ast}\left|\int x^\beta \curl u^{cc}(t,x) \,\dx\right|^p\right).\label{curl ut ucc}
\end{align}
For the first term on the right hand side of \eqref{curl ut ucc}, we use Lemma \ref{lemma: estimate time corrector} and the fact that $|M(t)|\approx \lambda^2$ and obtain
\begin{align*}
\sum_{x_0\in M(t)}\left(\frac{1}{\lambda\mu_1}\right)^2 \|\curl\tilde{\theta}_{x_0}(t,\cdot)\|^p_{L^\infty(\R^2)} &\leq C(R_0,u_0,e,\delta, \kappa,\varepsilon)\sum_{x_0\in M(t)}\left(\frac{1}{\lambda\mu_1}\right)^2\left(\frac{\lambda\mu_1\mu_2^2}{\omega}\right)^p\\
&\leq  C(R_0,u_0,e,\delta,\kappa,\varepsilon) \omega^{-p}\lambda^p\mu_1^{p-2}\mu_2^{2p}.
\end{align*}
For the second term in \eqref{curl ut ucc}, we use Lemma \ref{lemma: estimate correctors ii}, i.e.
\begin{align*}
\|\curl u^{cc}(t)\|^p_{L^\infty(\R^2)}&\leq C(R_0, u_0, e,\delta, \kappa,\varepsilon,N,M)\left(\omega^{-1}\lambda^{-M+1}\mu_1\mu_2^2 + \frac{1}{\omega}\right)^p.
\end{align*}
For the third term in \eqref{curl ut ucc}, we use integration by parts and Lemma \ref{lemma: estimate time corrector}
\begin{align*}
 \sum_{x_0\in M(t)} |\log(\lambda\mu_1)| \max_{|\beta|\leq N^\ast}\left|\int x^\beta\curl \tilde{\theta}_{x_0}(t,x) \,\dx\right|^p &= \sum_{x_0\in M(t)} |\log(\lambda\mu_1)| \max_{|\beta|\leq N^\ast}\left|\int -\nabla^\perp(x^\beta) \tilde{\theta}_{x_0}(t,x) \,\dx\right|^p\\
 &\leq C \sum_{x_0\in M(t)} |\log(\lambda\mu_1)|\|u^t\|^p_{L^1(\R^2)}\\
 &\leq  C(R_0,u_0,e,\delta, \kappa, \varepsilon)\lambda^2 |\log(\lambda\mu_1)| \omega^{-p}.
\end{align*}
Finally, for the last term in \eqref{curl ut ucc}, we use Lemma \ref{lemma: estimate correctors ii}
\begin{align*}
\max_{|\beta|\leq N^\ast}\left|\int x^\beta \curl u^{cc}(t,x) \,\dx\right|^p &\leq C \|\curl u^{cc}\|_{L^\infty(\R^2)}^p\\
&\leq  C(R_0, u_0, e,\delta, \kappa,\varepsilon,N,M)\left(\omega^{-1}\lambda^{-M+1}\mu_1\mu_2^2 + \frac{1}{\omega}\right)^p.
\end{align*}
Putting everything together yields the claim.
\end{proof}

\section{The new error}\label{sec: def of R}
This section is devoted to the definition of the new error $(r_1,R_1)$, which will be estimated in the next section.
\subsection{The new Reynolds-defect-equation}
Plugging  $u_1$ into the new Reynolds-defect-equation and writing 
\begin{align*}
u_1 = u_0 + w = u_0 + w^p +w^t = u_0 +u^p+u^t+u^c+u^{cc},
\end{align*} we need to define $(R_1, p_1)$ such that
\begin{align}
 - \dv \overset{\circ}{R_1} &= \dv(u_0\otimes (u_1 - u_0) + (u_1-u_0)\otimes u_0)\nonumber\\
&\hspace{0,3cm} +\dv(u^p\otimes (u_1-u_0-u^p)) + \dv((u_1-u_0 - u^p)\otimes u^p) \nonumber\\
&\hspace{0,3cm}  + \dv((u_1-u_0-u^p)\otimes (u_1-u_0-u^p))\nonumber\\
&\hspace{0,3cm} + \partial_t (u_1-u_0) + \dv(u^p\otimes u^p - \Rtr)+ \nabla(p_1 - p_0).\label{eq: iteration eq}
\end{align}
Let us analyse \eqref{eq: iteration eq}.
\subsection{\texorpdfstring{Analysis of the first three lines of \eqref{eq: iteration eq}}{Analysis of the first three lines of the iteration equation}}
Let us define
\begin{align*}
R^{\operatorname{lin},1} &= u_0\otimes (u_1 - u_0) + (u_1-u_0)\otimes u_0,\\
R^{\operatorname{lin},2}& =  u^p\otimes (u_1-u_0-u^p)+(u_1-u_0 - u^p)\otimes u^p,\\
R^{\operatorname{lin},3}&=  (u_1-u_0-u^p)\otimes (u_1-u_0-u^p),
\end{align*}
i.e.
\begin{align*}
&\dv(u_0\otimes (u_1 - u_0) + (u_1-u_0)\otimes u_0)\\
 &\hspace{0,3cm}+ \dv(u^p\otimes (u_1-u_0-u^p)) +\dv((u_1-u_0 - u^p)\otimes u^p) \\
&\hspace{0,3cm} + \dv((u_1-u_0-u^p)\otimes (u_1-u_0-u^p))\\
&= \dv(R^{\operatorname{lin},1} + R^{\operatorname{lin},2} + R^{\operatorname{lin},3}).
\end{align*}

\subsection{\texorpdfstring{Analysis of the fourth  line of \eqref{eq: iteration eq}}{Analysis of the fourth line of the iteration equation}}
Remembering Proposition \ref{prop: structure of perturbations} and going into the definition of the single perturbation parts, let us start by calculating
\begin{align}
\partial_t (u_1-u_0) + \dv(u^p\otimes u^p - \Rtr) +  \nabla(p_1 - p_0)&= \partial_t u^t + \partial_t u^{cc}_{I} + \dv(u^p\otimes u^p - \Rtr)\nonumber\\
&\hspace{0,3cm} + \partial_t(u^p + u^c)\nonumber\\
&\hspace{0,3cm} + \partial_t u^{cc}_{II}\nonumber\\
&\hspace{0,3cm} + \partial_t u^{cc}_{III}  +  \nabla(p_1 - p_0).\label{quadratic and time error}
\end{align}
We consider the first line on the right hand side of the previous calculation. By Proposition \ref{prop: building blocks} and the decomposition of the old error $\Rtr$ in \eqref{eq: R coeff}, we have
\begin{align}
\partial_t &u^t + \partial_t u^{cc}_{I} + \dv(u^p\otimes u^p - \Rtr)\nonumber\\
&=-\sum_k\partial_t a_k^2 \left(Y_k-\frac{1}{\omega}\xi_k\right)  -\sumk a_k^2 \partial_t Y_k\nonumber\\
&\hspace{0,3cm}+\dv\left(\sumk a_k^2 \left(W_k\otimes W_k - \frac{\xi_k}{|\xi_k|}\otimes \frac{\xi_k}{|\xi_k|}\right)\right) + \nabla(\chi_\kappa^2\rho)\nonumber\\
 &=-\sum_k\partial_t a_k^2 \left(Y_k-\frac{1}{\omega}\xi_k\right) +  \sumk a_k^2  \underbrace{\left(\dv\left(W_k\otimes W_k\right) - \partial_t Y_k\right)}_{=0}\nonumber\\
&\hspace{0,3cm} + \sumk  \left(W_k\otimes W_k-\frac{\xi_k}{|\xi_k|}\otimes \frac{\xi_k}{|\xi_k|}\right)\cdot\nabla a_k^2 + \nabla(\chi_\kappa^2\rho)\nonumber\\
&= \dv R^Y + \dv R^{\quadr} - \pi_1\nonumber\\
&\hspace{0,3cm}+\sum_k \left(W_k\otimes W_k-\frac{\xi_k}{|\xi_k|}\otimes \frac{\xi_k}{|\xi_k|}\right)\cdot\nabla a_k^2 -\partial_t a_k^2 \left(Y_k-\frac{1}{\omega}\xi_k\right)\nonumber\\
&\hspace{0,3cm}  - \sum_k\px\left(\left(W_k\otimes W_k-\frac{\xi_k}{|\xi_k|}\otimes \frac{\xi_k}{|\xi_k|}\right)\nabla a_k^2 -  \partial_t a_k^2 \left(Y_k-\frac{1}{\omega}\xi_k\right)\right)\label{quadr error i}
\end{align}
where we can directly define
\begin{align*}
R^{\quadr} &= \sum_k \tilde{\Rc}\left(\nabla a_k^2, W_k\otimes W_k - \frac{\xi_k}{|\xi_k|}\otimes \frac{\xi_k}{|\xi_k|}\right),\\
R^Y &=  -\sum_k \Rc\left(\partial_t a_k^2, Y_k-\frac{1}{\omega}\xi_k\right),\\
\pi_1 &=  -\nabla (\chi_\kappa^2\rho).
\end{align*}
Note that this is well-defined since $\int_{\tor} W_k\otimes  W_k (t,x)\,\dx = \frac{\xi_k}{|\xi_k|}\otimes\frac{\xi_k}{|\xi_k|} $ and  $\int_{\tor} Y_k (t,x)\,\dx =\frac{1}{\omega}\xi_k$  by Proposition \ref{prop: building blocks}.
For the second line in \eqref{quadratic and time error}, we take into account Proposition \ref{prop: structure of perturbations} and the elementary identity $\nabla^\perp\curl = \Delta - \nabla\dv$ and  obtain
\begin{align}
\partial_t(u^p + u^c) &= \sumk \nabla^\perp\curl\Delta^N\dv \partial_t \left(a_k (A_k+A_k^T)\right)\nonumber\\
&= \sum_k \dv\left(\Delta^{N+1} \partial_t\left(a_k ( A_k +  A_k^T)\right)\right) - \nabla\dv \Delta^N\dv\partial_t\left(a_k(A_k+A_k^T)\right)\nonumber\\
&= \dv R^{\tme,1} - \pi_2\label{time error i}
\end{align}
with
\begin{align*}
R^{\tme,1} &= \sum_k \Delta^{N+1} \partial_t\left(a_k (A_k + A_k^T)\right),\\
\pi_2 &= \nabla\sum_k\dv \Delta^N\dv\partial_t\left(a_k(A_k+A_k^T)\right).
\end{align*}
The third line in \eqref{quadratic and time error} can be written as
\begin{align*}
 \partial_t u^{cc}_{II} = \dv R^{\tme,2} +  \partial_t u^{cc}_{II} - \px( \partial_t u^{cc}_{II})
\end{align*}
where
\begin{align*}
R^{\tme,2} = \antidvx\left(\px ( \partial_t u^{cc}_{II})\right)
\end{align*}
and $\antidvx$ is the operator from Proposition \ref{thm: isett}. In the fifth line of \eqref{quadratic and time error}, we define
\begin{align}\label{time error ii}
\pi_3= -\partial_t u^{cc}_{III}.
\end{align}
Collecting \eqref{quadr error i}, \eqref{time error i} and \eqref{time error ii}, we obtain from \eqref{quadratic and time error}
\begin{align}
\partial_t &(u_1-u_0) + \dv(u^p\otimes u^p - \Rtr) +  \nabla(p_1 - p_0)\nonumber\\
&=\dv R^Y + \dv R^{\quadr}  + \dv R^{\tme,1} + \dv R^{\tme,2}\nonumber\\
&\hspace{0,3cm}+ \nabla(p_1-p_0) -\pi_1-\pi_2 -\pi_3\nonumber\\
&\hspace{0,3cm}+\sum_k \left(W_k\otimes W_k-\frac{\xi_k}{|\xi_k|}\otimes \frac{\xi_k}{|\xi_k|}\right)\cdot\nabla a_k^2 -\partial_t a_k^2 \left(Y_k-\frac{1}{\omega}\xi_k\right) + \partial_tu^{cc}_{II}\nonumber\\
&\hspace{0,3cm}  - \sum_k\px\left(\left(W_k\otimes W_k-\frac{\xi_k}{|\xi_k|}\otimes \frac{\xi_k}{|\xi_k|}\right)\cdot\nabla a_k^2 -  \partial_t a_k^2 \left(Y_k-\frac{1}{\omega}\xi_k\right)+ \partial_tu^{cc}_{II}\right).\label{quadr and time error ii}
\end{align}
We note that $-\pi_1-\pi_2-\pi_3$ is a compactly supported gradient. Taking the integral \eqref{quadr and time error ii} in the previous equality and using that our perturbations are higher derivatives of compactly supported functions, we get
\begin{align*}
0 = \int \sum_k \left(W_k\otimes W_k-\frac{\xi_k}{|\xi_k|}\otimes \frac{\xi_k}{|\xi_k|}\right)\cdot\nabla a_k^2 -\partial_t a_k^2 \left(Y_k-\frac{1}{\omega}\xi_k\right) + \partial_tu^{cc}_{II} \,\dx.
\end{align*}
In the same way, integrating over the angular momentum in \eqref{quadr and time error ii} and using that for any compactly supported matrix as well as for any compactly supported gradient, the integral of the angular momentum vanishes, we obtain
\begin{align*}
0 &= \int x_1 \left(\sum_k \left(W_k\otimes W_k-\frac{\xi_k}{|\xi_k|}\otimes \frac{\xi_k}{|\xi_k|}\right)\cdot\nabla a_k^2 -\partial_t a_k^2 \left(Y_k-\frac{1}{\omega}\xi_k\right) + \partial_tu^{cc}_{II}\right)_2 \,\dx\\
&- \int x_2 \left(\sum_k \left(W_k\otimes W_k-\frac{\xi_k}{|\xi_k|}\otimes \frac{\xi_k}{|\xi_k|}\right)\cdot\nabla a_k^2 -\partial_t a_k^2 \left(Y_k-\frac{1}{\omega}\xi_k\right) + \partial_tu^{cc}_{II}\right)_1 \,\dx.
\end{align*}
Therefore, 
\begin{align*}
\sum_k \left(W_k\otimes W_k-\frac{\xi_k}{|\xi_k|}\otimes \frac{\xi_k}{|\xi_k|}\right)\cdot\nabla a_k^2 -\partial_t a_k^2 \left(Y_k-\frac{1}{\omega}\xi_k\right) + \partial_tu^{cc}_{II}\in X
\end{align*}
with the notation from Proposition \ref{thm: isett}, where $X$ is the space of smooth vector fields supported in $B_1(0)$ with vanishing linear and angular momenta. Consequently, \eqref{quadr and time error ii} reduces to
\begin{align*}
\partial_t (u_1-u_0) + \dv(u^p\otimes u^p - \Rtr) +  \nabla(p_1 - p_0)&=\dv R^Y + \dv R^{\quadr}  + \dv R^{\tme,1} + \dv R^{\tme,2}\nonumber\\
&\hspace{0,3cm}+ \nabla(p_1-p_0) -\pi_1-\pi_2 -\pi_3.
\end{align*}

\subsection{Definition of the new error and the new pressure}
Altogether, we define
\begin{align*}
R_1 &= -\left(R^{\operatorname{lin},1} + R^{\operatorname{lin},2} +R^{\operatorname{lin},3}  +R^{\quadr}+ R^Y  + R^{\tme,1} + R^{\tme,2} \right),\\
\nabla p_1 &= \nabla p_0 + \pi_1 + \pi_2+ \pi_3 +\frac{1}{2}\nabla\operatorname{tr}R_1.
\end{align*}
Recall that, by definition, $\pi_1,\pi_2,\pi_3$ are gradients.
With that definition, it holds
\begin{align*}
-\dv  \overset{\circ}{R_1} = \partial_t u_1 + \dv(u_1 \otimes u_1) + \nabla p_1.
\end{align*}

\section{Estimates of the new error}\label{sec: error estimates}
We will now estimate the different parts of $R_1$ that were defined in the previous section.
\begin{lemma}[Estimate of $R^{\operatorname{lin},1} $]\label{lemma: est rlin}
It holds
\begin{align*}
\|R^{\operatorname{lin},1}(t) \|_{L^1(\R^2)}\leq  C(R_0,u_0,e,\delta,\varepsilon,N,M)\left(\frac{\mu_1}{\mu_2} + \frac{\mu_1^{\frac{1}{2}}\mu_2^{\frac{1}{2}}}{\omega}+ \mu_1^{-\frac{1}{2}}\mu_2^{-\frac{1}{2}}\right)
\end{align*}
uniformly in $t$.
\end{lemma}

\begin{proof}
Using Hölder's inequality and Lemma \ref{lemma: estimate up}, Lemma \ref{lemma: estimate correctors i},  Lemma \ref{lemma: estimate time corrector} and Lemma \ref{lemma: estimate correctors ii}, we have
\begin{align*}
\|R^{\operatorname{lin},1}(t)\|_{L^1(\R^2)} &\leq 2 \|u_0(t)\|_{L^2(\R^2)}\left( \|u^t(t)\|_{L^2(\R^2)}+ \|u^c(t)\|_{L^2(\R^2)}+\|u^{cc}(t)\|_{L^2(\R^2)} \right)\\
&\hspace{0,3cm}   + 2  \|u_0(t)\|_{L^\infty(\R^2)}\|u^p(t)\|_{L^1(\R^2)} \\
&\leq C(R_0,u_0,e,\delta,\varepsilon,N,M)\left(\frac{\mu_1}{\mu_2} + \frac{\mu_1^{\frac{1}{2}}\mu_2^{\frac{1}{2}}}{\omega}+ \mu_1^{-\frac{1}{2}}\mu_2^{-\frac{1}{2}}\right) .\qedhere
\end{align*}
\end{proof}

\begin{lemma}[Estimate of $R^{\operatorname{lin},2} $]\label{lemma: est rlin2}
It holds
\begin{align*}
\|R^{\operatorname{lin},2}(t) \|_{L^1(\R^2)}\leq  C(R_0,r_0,e,\delta,\varepsilon,N,M) \left(\frac{\mu_1}{\mu_2} + \frac{\mu_1^{\frac{1}{2}}\mu_2^{\frac{1}{2}}}{\omega}\right)
\end{align*}
uniformly in $t$.
\end{lemma}

\begin{proof}
We have by  Lemma \ref{lemma: estimate up}, Lemma \ref{lemma: estimate correctors i},  Lemma \ref{lemma: estimate time corrector} and Lemma \ref{lemma: estimate correctors ii}
\begin{align*}
\|R^{\operatorname{lin},2}(t) \|_{L^1(\R^2)}&\leq 2 \|u^p(t)\|_{L^2(\R^2)}\left(  \|u^t(t)\|_{L^2(\R^2)} + \|u^c(t)\|_{L^2(\R^2)}  +\|u^{cc}(t)\|_{L^2(\R^2)}\right)\\
&\leq C(R_0, u_0, e,\delta, \kappa,\varepsilon,N,M) \left(\frac{\mu_1}{\mu_2} + \frac{\mu_1^{\frac{1}{2}}\mu_2^{\frac{1}{2}}}{\omega}\right).
\end{align*}
\end{proof}

\begin{lemma}[Estimate of $R^{\operatorname{lin},3} $]\label{lemma: est rlin3}
It holds
\begin{align*}
\|R^{\operatorname{lin},3}(t) \|_{L^1(\R^2)}\leq  C(R_0,u_0,e,\delta,\varepsilon,N,M)\left(\left(\frac{\mu_1}{\mu_2}\right)^2 + \left(\frac{\mu_1^\frac{1}{2}\mu_2^\frac{1}{2}}{\omega}\right)^2\right)
\end{align*}
uniformly in $t$.
\end{lemma}

\begin{proof}
Since $$R^{\operatorname{lin},3}=  (u_1-u_0-u^p)\otimes (u_1-u_0-u^p) = ( u^t + u^c + u^{cc})\otimes ( u^t + u^c + u^{cc}),$$ we have
\begin{align*}
\|R^{\operatorname{lin},3}(t) \|_{L^1(\R^2)}\leq 4\left(\|u^t(t)\|_{L^2(\R^2)}^2  + \|u^c(t)\|^2_{L^2(\R^2)} + \|u^{cc}(t)\|^2_{L^2(\R^2)}\right),
\end{align*}
hence the assertion follows from Lemma \ref{lemma: estimate correctors i},  Lemma \ref{lemma: estimate time corrector} and Lemma \ref{lemma: estimate correctors ii}.
\end{proof}

\begin{lemma}[Estimate of $R^{\quadr}$]\label{lemma: est rquadr}
It holds 
\begin{align*}
\|R^{\quadr}(t)\|_{L^1(\R^2)}\leq \frac{C(R_0,u_0,e,\delta,\varepsilon)}{\lambda}
\end{align*}
uniformly in $t$.
\end{lemma}

\begin{proof}
We remember that $R^{\quadr} = \sum_k \tilde{\Rc}\left(\nabla a_k^2, W_k\otimes W_k - \frac{\xi_k}{|\xi_k|}\otimes \frac{\xi_k}{|\xi_k|}\right)$.
The estimate follows directly from Remark \ref{rem: improved antidiv}, the scaling of $W_k$ (see Proposition \ref{prop: building block estimates}) and the estimates on $a_k$ in \eqref{est ak}.
\end{proof}

\begin{lemma}[Estimate of $R^{Y}$]\label{lemma: est ry}
It holds 
\begin{align*}
\|R^{Y}(t)\|_{L^1(\R^2)}\leq \frac{C(R_0,u_0,e,\delta,\varepsilon)}{\omega\lambda}
\end{align*}
uniformly in $t$.
\end{lemma}
\begin{proof}
As for $R^{\quadr}$, this is a direct application of Remark \ref{rem: improved antidiv} and Proposition \ref{prop: building block estimates}.
\end{proof}

\begin{lemma}[Estimate of $R^{\tme,1}$]\label{lemma: est rtime}
It holds
\begin{align*}
\|R^{\tme,1}(t)\|_{L^1(\R^2)} \leq C(R_0,u_0,e,\delta,\varepsilon) \omega \mu_1^\frac{1}{2}\mu_2^{-\frac{3}{2}}
\end{align*}
uniformly in $t$.
\end{lemma}

\begin{proof}
$R^{\tme,1}$ was defined as $R^{\tme,1} = \sum_k \Delta^{N+1} \partial_t\left(a_k (A_k + A_k^T)\right)$. By Proposition \ref{prop: building block estimates} and \eqref{est ak} this yields
\begin{align*}
\|R^{\tme,1}(t)\|_{L^1(\R^2)}& \leq C(R_0,u_0,e,\delta,\varepsilon) \sum_k \|\partial_t A_k\|_{W^{2N+2,1}(\tor)}\\
 &\leq C(R_0,u_0,e,\delta,\varepsilon) \omega\mu_1^{\frac{1}{2}}\mu_2^{-\frac{3}{2}}.
\end{align*}
\end{proof}

\begin{lemma}[Estimate of $R^{\tme,2}$]\label{lemma: est rtime 2}
It holds
\begin{align*}
\|R^{\tme,2}(t)\|_{L^1(\R^2)} \leq C(R_0, u_0, e,\delta, \kappa,\varepsilon,N,M)\left(\lambda^{-M+1}\mu_1 + \frac{1}{\omega}\right).
\end{align*}
uniformly in $t$.
\end{lemma}

\begin{proof}
It holds $R^{\tme,2} = \antidvx\left(\px ( \partial_t u^{cc}_{II})\right)$. By Proposition \ref{thm: isett}, \eqref{est: px} and Lemma \ref{lemma: estimate correctors ii}
\begin{align*}
\|R^{\tme,2}(t)\|_{L^1(\R^2)} \leq C \|\partial_t u^{cc}_{II}\|_{L^1(\R^2)} \leq C(R_0, u_0, e,\delta, \kappa,\varepsilon,N,M)\left(\lambda^{-M+1}\mu_1 + \frac{1}{\omega}\right).
\end{align*}
\end{proof}

\section{Proof of the main proposition}\label{sec: proof of main prop}
Proposition \ref{prop: main proposition} is proved by choosing all the parameters appropriately, which we do in this section. Let us set 
\begin{itemize}
\item $\mu_1= \lambda^\alpha$,
\item $\mu_2 = \lambda \mu_1= \lambda^{1+\alpha}$,
\item $\omega = \lambda^\beta$
\end{itemize}
for some $\alpha,\beta>0$ to be chosen below. We collect the estimates from Section \ref{sec: estimates perturbations}, \ref{sec: curl estimates} and \ref{sec: error estimates}  where the parameters $\mu_1,\mu_2$ and $\omega$ need to be balanced in Table \ref{tab: estimates}.
\begin{table}[h]
\caption{All quantities that have to be controlled in terms of oscillation and concentration parameters and the phase speed}
\begin{tabular}[h]{l|l|l|l}
Lemma & Term & Order & = \\
\hline
\ref{lemma: estimate time corrector}, \ref{lemma: estimate correctors ii}, \ref{lemma: energy inc},  & $u^t, u^{cc}$ in $L^2(\R^2)$, Energy increment, & $\omega^{-1}\mu_1^\frac{1}{2}\mu_2^\frac{1}{2}$ & $\lambda^{-\beta +\alpha + \frac{1}{2}}$\\
\ref{lemma: est rlin}--\ref{lemma: est rlin3} & $R^{\operatorname{lin},1}$, $R^{\operatorname{lin},2}$, $R^{\operatorname{lin},3}$ in $L^1(\R^2)$ & &\\
\ref{lemma: estimate correctors i}, \ref{lemma: energy inc}, & $u^c$ in $L^2(\R^2)$, Energy increment & $\mu_1\mu_2^{-1}$ & $\lambda^{-1}$\\
\ref{lemma: est rlin}--\ref{lemma: est rlin3} & $R^{\operatorname{lin},1}$, $R^{\operatorname{lin},2}$, $R^{\operatorname{lin},3}$ in $L^1(\R^2)$ & &\\
\ref{lemma: curl up} & $\curl (w^p)$ in $H^p(\R^2)$ & $\lambda\mu_1^{\frac{1}{2}-\frac{2}{p}}\mu_2^\frac{3}{2}$ & $\lambda^{\frac{5}{2} + \alpha(2-\frac{2}{p})}$\\

\ref{lemma: curl ut} & $\curl (w^t)$ in $H^p(\R^2)$ & $\omega^{-1}\lambda\mu_1^{1-\frac{2}{p}}\mu_2^2$, & $\lambda^{3-\beta +\alpha(3-\frac{2}{p})}$\\
&& $\omega^{-1}\lambda^{-M+1}\mu_1\mu_2^2$,& $\lambda^{-\beta+3\alpha+3-M}$\\
&& $\lambda^2\log(\lambda\mu_1)\omega^{-p}$ & $(\alpha+1)\log(\lambda)\lambda^{2-\beta p}$\\
\ref{lemma: est rtime} & $R^{\tme,1}$ in $L^1(\R^2)$ & $\omega\mu_1^\frac{1}{2}\mu_2^{-\frac{3}{2}}$ & $\lambda^{\beta-\alpha-\frac{3}{2}}$\\
\ref{lemma: est rtime 2} & $R^{\tme,2}$ in $L^1(\R^2)$ & $\lambda^{-M+1}\mu_1$ & $\lambda^{-M+1+\alpha}$
\end{tabular}\label{tab: estimates}
\end{table}
\vspace{0,3cm}  \\We choose $\alpha,\beta$ and $M$ such that all the exponents in the fourth column of the previous table are negative. This is clear for the second row. Since $2-\frac{2}{p}<0$, we can choose $\alpha \gg 1$ so large such that
\begin{align*}
\frac{5}{2} + \alpha(2-\frac{2}{p}) &< 0.
\end{align*}
This yields negative exponents in Line 3. Furthermore, since $3-\frac{2}{p}<1$, let us choose $\alpha$ large enough such that 
\begin{align*}
3+ \alpha(3-\frac{2}{p})&<\alpha+\frac{1}{2},\\
\frac{3}{p}&< \alpha + \frac{1}{2}.
\end{align*}
With this choice of $\alpha$, we only need $\beta$ to satisfy
\begin{align*}
3+\alpha(3-\frac{2}{p})<\alpha + \frac{1}{2} <\beta < \alpha +\frac{3}{2}.
\end{align*}
With such a $\beta$, Line 1 -- 4  and Line 7 in Table \ref{tab: estimates} have negative exponents of $\lambda$. Also, this yields in Line 6
\begin{align*}
(\alpha+1)\log(\lambda)\lambda^{2-\beta p} \leq (\alpha+1)\lambda^{3-\beta p}
\end{align*}
has a negative exponent since $\beta >\frac{3}{p}$. Having $\alpha$ and $\beta$ fixed, it only remains to choose $M$. Since $M$ enters all the remaining exponents with a negative sign, we can simply pick $M\in\N$ large enough such that all exponents are negative.
Let 
\begin{align*}
\gamma_0 = \text{ exponent in the table with the smallest magnitude}<0
\end{align*}
which is negative by our choice of $\alpha,\beta, M$.
We can now verify the claims of Proposition \ref{prop: main proposition}. For $(i)$, we have by \eqref{energy increment}
\begin{align*}
\left|e(t) \left(1-\frac{\delta}{2}\right) - \int_{\R^2} |u_1|^2\,\dx\right|&<  \frac{1}{8}\delta e(t) +  C(R_0,r_0, u_0,e, \delta,\varepsilon)\left(\lambda^{\gamma_0} + \mu_1^{-\frac{1}{2}}\mu_2^{-\frac{1}{2}}\right).
\end{align*}
and we can choose $\lambda$ large enough such that $(i)$ is satisfied. For $(iii)$, we use Lemma \ref{lemma: estimate up}, \eqref{est up with eps},  Lemma \ref{lemma: estimate time corrector}, Lemma \ref{lemma: estimate correctors i} and Lemma \ref{lemma: estimate correctors ii}
\begin{align*}
\|u_1(t)-u_0(t)\|_{L^2(\R^2)}&\leq \|u^p(t)\|_{L^2(\R^2)} +  \|u^t(t)\|_{L^2(\R^2)}  +\|u^c(t)\|_{L^2(\R^2)}  +\|u^{cc}(t)\|_{L^2(\R^2)}\\
&\leq 40\delta^\frac{1}{2} +\frac{ C(R_0,u_0,e,\delta,\varepsilon)}{\lambda^\frac{1}{2}} + C(R_0,u_0,e,\delta,\varepsilon)\lambda^{\gamma_0}.
\end{align*}
We can choose $\lambda$ large enough such that 
\begin{align*}
\|u_1-u_0(t)\|_{L^2(\R^2)}&\leq 41 \delta^\frac{1}{2} ,
\end{align*}
i.e. $(v)$ is satisfied with $M_0 = 41$.
For $(iv)$, we use Lemmas \ref{lemma: curl up} and \ref{lemma: curl ut}:
\begin{align*}
\|\curl (u_1 - u_0)(t)\|^p_{H^p(\R^2)} &\leq \|\curl (w^p) (t)\|^p_{H^p(\R^2)} + \|\curl (w^t)(t)\|^p_{H^p(\R^2)} \\
&\leq C(R_0,u_0,e,\delta,\kappa,\varepsilon)\lambda^{p\gamma_0} 
\end{align*}
and we can choose $\lambda$ large enough such that $(iv)$ is satisfied. For $(ii)$, we have by Lemma \ref{lemma: est rlin},  \ref{lemma: est rlin2},  \ref{lemma: est rlin3},  \ref{lemma: est rquadr}, \ref{lemma: est ry}, \ref{lemma: est rtime} and \ref{lemma: est rtime 2}
\begin{align*}
\|R_1(t)\|_{L^1(\R^2)}&\leq \|R^{\operatorname{lin},1}(t)\|_{L^1(\R^2)} + \|R^{\operatorname{lin},2}(t)\|_{L^1(\R^2)} + \|R^{\operatorname{lin},3}(t)\|_{L^1(\R^2)}
\\
 &\hspace{0,3cm} +  \|R^{\quadr}(t)\|_{L^1(\R^2)} + \|R^Y(t)\|_{L^1(\R^2)}  + \|R^{\tme,1}(t)\|_{L^1(\R^2)} + \|R^{\tme,2}(t)\|_{L^1(\R^2)}\\
&\leq C(R_0,u_0,e,\delta,\kappa,\varepsilon,N) \lambda^{\gamma_0},
\end{align*}
hence we can choose $\lambda$ large enough to obtain $(ii)$. This proves $(i)$--$(iv)$ in Proposition \ref{prop: main proposition}. Assume we are given two energy profiles $e_1,e_2$ with $e_1=e_2$ on $[0,t_0]$ for some $t_0\in[0,1]$.  The values that we add with $w(t)$ depend only on pointwise (in time) values of the previous steps. Therefore, one can do the construction for $e_1$ and $e_2$ simultaneously, choosing the same values for all the parameters in each iteration step, thereby producing two solutions $u_1$, $u_2$ to \eqref{2D Euler} that satisfy $u_1=u_2$ on $[0,t_0]$. This concludes the proof of Proposition \ref{prop: main proposition} and therefore also for Theorem \ref{thm:main}.

\bibliographystyle{abbrvdin}
\bibliography{bibliography}
\end{document}